\documentclass[11pt,reqno,tbtags]{amsart}
\usepackage{amssymb}
\usepackage[utf8]{inputenc} 
\usepackage[T1]{fontenc}  
\usepackage{url}
\usepackage[square,numbers]{natbib}

\numberwithin{equation}{section}

\allowdisplaybreaks

\newtheorem{theorem}{Theorem}[section]
\newtheorem{lemma}[theorem]{Lemma}

\newtheorem{corollary}[theorem]{Corollary}

\theoremstyle{definition}
\newtheorem{example}[theorem]{Example}
\newtheorem{definition}[theorem]{Definition}
\newtheorem{problem}[theorem]{Problem}
\newtheorem{remark}[theorem]{Remark}

\newtheorem*{ack}{Acknowledgement}

\newtheorem{divisor}{Divisor method, formulation}

\newtheorem{kvot}{Quota method, formulation}

\theoremstyle{remark}

{\begin{math}}
 {\end{math}}

\newenvironment{romenumerate}[1][0pt]{
\addtolength{\leftmargini}{#1}\begin{enumerate}
 }{\end{enumerate}}

\newcounter{oldenumi}
{\setcounter{oldenumi}{\value{enumi}}
\begin{romenumerate} \setcounter{enumi}{\value{oldenumi}}}
{\end{romenumerate}}

\newcounter{thmenumerate}

\newcounter{romxenumerate}   

\newcounter{xenumerate}   

\newcommand{\refT}[1]{Theorem~\ref{#1}}
\newcommand{\refC}[1]{Corollary~\ref{#1}}
\newcommand{\refL}[1]{Lemma~\ref{#1}}
\newcommand{\refR}[1]{Remark~\ref{#1}}
\newcommand{\refS}[1]{Section~\ref{#1}}
\newcommand{\refP}[1]{Problem~\ref{#1}}
\newcommand{\refD}[1]{Definition~\ref{#1}}
\newcommand{\refE}[1]{Example~\ref{#1}}

\newcommand{\refApp}[1]{Appendix~\ref{#1}}

\newcommand{\refTab}[1]{Table~\ref{#1}}
\newcommand{\refand}[2]{\ref{#1} and~\ref{#2}}
\newcommand\refform[1]{formulation \ref{#1}}
\newcommand\refForm[1]{Formulation \ref{#1}}

\begingroup
  \divide\count255 by 60
  \multiply\count255 by -60
  \advance\count255 by \time
  \ifnum \count255 < 10 \xdef\klockan{\the\count1.0\the\count255}
  \else\xdef\klockan{\the\count1.\the\count255}\fi
\endgroup

\newcommand\nopf{\qed}   

\newcommand{\sumim}{\sum_{i=1}^m}

\newcommand{\sumjk}{\sum_{j=1}^k}

\newcommand{\sumjm}{\sum_{j=1}^m}
\newcommand{\sumjjm}{\sum_{j=2}^m}
\newcommand{\sumjmi}{\sum_{j=1}^{m-1}}
\newcommand{\sumkm}{\sum_{k=1}^m}

\newcommand\set[1]{\ensuremath{\{#1\}}}

\newcommand\xpar[1]{(#1)}
\newcommand\bigpar[1]{\bigl(#1\bigr)}
\newcommand\Bigpar[1]{\Bigl(#1\Bigr)}
\newcommand\biggpar[1]{\biggl(#1\biggr)}
\newcommand\lrpar[1]{\left(#1\right)}

\newcommand\lrsqpar[1]{\left[#1\right]}
\newcommand\xcpar[1]{\{#1\}}

\newcommand\abs[1]{|#1|}
\newcommand\bigabs[1]{\bigl|#1\bigr|}

\newcommand\biggabs[1]{\biggl|#1\biggr|}
\newcommand\Biggabs[1]{\Biggl|#1\Biggr|}
\newcommand\lrabs[1]{\left|#1\right|}
\def\rompar(#1){\textup(#1\textup)}    
\newcommand\xfrac[2]{#1/#2}

\def\xexp(#1){e^{#1}}
\newcommand\ceil[1]{\lceil#1\rceil}
\newcommand\floor[1]{\lfloor#1\rfloor}
\newcommand\lrfloor[1]{\left\lfloor#1\right\rfloor}
\newcommand\frax[1]{\{#1\}}
\newcommand\Bigfrax[1]{\Bigl\{#1\Bigr\}}
\newcommand\lrfrax[1]{\left\{#1\right\}}

\newcommand\ntoo{\ensuremath{{n\to\infty}}}
\newcommand\Ntoo{\ensuremath{{N\to\infty}}}
\newcommand\mtoo{\ensuremath{{m\to\infty}}}

\newcommand\ktoo{\ensuremath{{k\to\infty}}}
\newcommand\xtoo{\ensuremath{{x\to\infty}}}
\newcommand\too{\ensuremath{{\to\infty}}}

\newcommand\punkt[1]{\if.#1\else.\spacefactor1000\fi{#1}}
\newcommand\iid{i.i.d\punkt}    
\newcommand\ie{i.e\punkt}
\newcommand\eg{e.g\punkt}

\newcommand\cf{cf\punkt}

\newcommand{\etc}{etc\punkt}

\newcommand\ii{\mathrm{i}}

\newcommand{\tend}{\longrightarrow}
\newcommand\dto{\overset{\mathrm{d}}{\tend}}
\newcommand\pto{\overset{\mathrm{p}}{\tend}}
\newcommand\ptox{\overset{\mathrm{p*}}{\tend}}

\newcommand\eqd{\overset{\mathrm{d}}{=}}

\newcommand\bbR{\mathbb R}

\newcommand\bbN{\mathbb N}
\newcommand\bbT{\mathbb T}
\newcommand\bbQ{\mathbb Q}
\newcommand\bbZ{\mathbb Z}

\newcounter{CC}
\newcounter{cc}

\newcommand\E{\operatorname{\mathbb E{}}}
\renewcommand\P{\operatorname{\mathbb P{}}}
\newcommand\Var{\operatorname{Var}}
\newcommand\Cov{\operatorname{Cov}}

\newcommand\Exp{\operatorname{Exp}}

\newcommand\sgn{\operatorname{sgn}}

\newcommand\ga{\alpha}
\newcommand\gb{\beta}
\newcommand\gd{\delta}
\newcommand\gD{\Delta}
\newcommand\gf{\varphi}
\newcommand\gam{\gamma}
\newcommand\gG{\Gamma}

\newcommand\gl{\lambda}

\newcommand\eps{\varepsilon}

\newcommand\ett[1]{\boldsymbol1\xcpar{#1}}

\newcommand\qw{^{-1}}
\newcommand\qww{^{-2}}
\newcommand\qq{^{1/2}}

\newcommand\qqc{^{3/2}}

\renewcommand{\=}{:=}

\newcommand\intoi{\int_0^1}
\newcommand\intoo{\int_0^\infty}

\newcommand\oi{[0,1]}
\newcommand\oix{[0,1)}

\newcommand\dd{\,\textup{d}}

\newcommand\ppm{\ensuremath{p_1,\dots,p_m}}
\newcommand\linQ{linearly independent over\/ $\bbQ$}
\newcommand\StL{Sainte-Laguë}
\newcommand\rd[1]{\left[#1\right]_d}
\newcommand\rga[1]{\left[#1\right]_\ga}
\newcommand\rgb[1]{\left[#1\right]_\gb}
\newcommand\rgx[2]{\left[#2\right]_{#1}}
\newcommand\U{\mathrm{U}}
\newcommand\Modm{\operatorname{Mod}_m}
\newcommand\tell{\tilde\ell}
\newcommand\gamq{$\gam$-quota}
\newcommand\gamqm{$\gam$-quota method}
\newcommand\tU{\widetilde U}
\newcommand\gbdiv{$\gb$-linear divisor method}
\newcommand\gar{$\ga$-rounding}
\newcommand\empho{\relax}
\newcommand\Lx{\bar L}
\newcommand\bX{\bar X}
\newcommand\hX{\widehat X}
\newcommand\tX{\widetilde X}
\newcommand\bY{\bar Y}
\newcommand\hY{\widehat Y}
\newcommand\ox[1]{_{(#1)}}
\newcommand\oy[1]{_{[#1]}}
\newcommand\sxm{\mathfrak S_m}
\newcommand\sxmge{\mathfrak S_{m,\ge}}
\newcommand\sxx[1]{\mathfrak S_{#1}}
\newcommand\sld{\tilde S}
\newcommand\sldlim{S^*}
\newcommand\nn{n_0}
\newcommand\nntoo{\ensuremath{{\nn\to\infty}}}
\newcommand\NN{N_0}
\newcommand\NNtoo{\ensuremath{{\NN\to\infty}}}
\newcommand\Ntoox{\ensuremath{{N\ptox\infty}}}
\newcommand\sij{s_{ij}}
\newcommand\gdij{\gD_{ij}}
\newcommand\sbm{S_{\gb,m}}
\newcommand\sgm{\bar S_{\gam,m}}
\newcommand\gbu{\hat U}
\newcommand\XX{Z}
\newcommand\tjr{t_{jr}}
\newcommand\maxim{\max_{1\le i\le m}}
\newcommand\minim{\min_{1\le i\le m}}
\newcommand\hh{^{\mathrm H}}
\newcommand\jj{^{\mathrm J}}
\newcommand\qa{^{\mathrm A}}
\newcommand\sjj{S\jj}
\newcommand\saa{|S\qa|}
\newcommand\sjjm{S_m\jj}
\newcommand\saam{|S_m\qa|}
\newcommand\yy{{\mathbf y}}

\newcommand\app{apparentement}

\newcommand\REM[1]{{\raggedright\texttt{[#1]}\par\marginal{XXX}}}

\newcommand\citetq[2]{\citeauthor{#2} \cite[{\frenchspacing #1}]{#2}}

\newcommand\urladdrx[1]{{\urladdr{\def~{{\tiny$\sim$}}#1}}}

\begin{document}
\title
{Asymptotic bias of some election methods}

\date{28 October 2011} 

\author{Svante Janson}
\address{Department of Mathematics, Uppsala University, PO Box 480,
SE-751~06 Uppsala, Sweden}
\email{svante.janson@math.uu.se}
\urladdrx{http://www.math.uu.se/~svante/}


\begin{abstract}
Consider an election where $N$ seats are distributed among parties with 
proportions $\ppm$ of the votes.
We study, for the common divisor and quota methods, 
the asymptotic distribution, and in particular the mean, of the
seat excess of a party, \ie{} the
difference between the number of seats given to the party and the 
(real) number $Np_i$ that yields exact proportionality.
Our approach is to keep $\ppm$ fixed and let $N\to\infty$, with $N$ random
in a suitable way.

In particular, we give formulas showing the bias favouring large or small
parties for the different election methods.
\end{abstract}

\maketitle

\section{Introduction}\label{S:intro}

The idea of a proportional election method is that each party gets a
number of seats that is proportional to the number of votes.
The same mathematical problem arises if seats are to be apportioned
(before the election) between states in a union or between
multi-member constituences according to their populations.
In particular, the problem of apportionment to the House of Representatives
in the United States has been the source of much debate as well as much
research since 1790, 
including the invention of several  important election methods;
see \citet{BY} for a detailed history.
(The European Union has yet to agree on a similar procedure;
see \citet{CambridgeCompromise} for a recent
attempt to replace the current 
political bartering with a ``mathematical formula'', \ie{} a well-defined
method to allocate seats in the European parliament to the member
states.) 
To fix the terminology, we talk in this paper about parties in an election;
the case of apportionment between states \etc{} differs mathematically only in
language. 

Of course, exact proportionality  is in general not possible, since the
number of seats has to be an integer, and the election method can be seen as
a method to ``round'' the real numbers  that yield exact
proportionality to integers. 
The purpose of this paper is to study the resulting deviation from perfect
proportionality, \ie{} the ``rounding errors'' introduced by
the method. In particular, we are interested in whether there is a systematic
bias favouring  large or small parties for various election methods.

This problem has received a considerable amount of attention going back to
(at least) 
\citet{Polya1918,Polya1918b,Polya1919}. 
(\citet{StL2} did  a related comparison between two different methods
already in 1910, see \refE{EStL}.)
Some
investigations have analysed data from real elections or from simulations
(with party sizes chosen at random according to some predetermined
distributions), see \eg{} \citet{BY} and \citet{Pukel2003b}.
The problem has also been studied theoretically by assuming
that the parties have random sizes with a uniform distribution of the
relative sizes over all possibilities, see 
\citet{Polya1918,Polya1918b,Polya1919},
\citet{Pukel2003b},
\citet{HPS2005}
and 
Drton and Schwingenschlögl \cite{DSbias2004,DSbias2005,DSvariance}.
Random party sizes with more general distributions have been considered by
\citet{HPS2004,HPS2005} and
\citet{Sch2008}.
(See also \refS{Srandom} below.)

The approach in the present paper is different.
We consider a number of parties with fixed sizes
and let the total number of seats $N$ be random. More precisely, we let $N$
be chosen uniformly at random in $1,2,\dots,\NN$ for some large integer $\NN$,
and then we take the limit as $\NN\too$. (Limits obtained in this way, if they
exist, can be regarded as results for ``a random positive integer'' $N$. See
also \refR{RN}.)
The same approach has been used in \cite{SJ258} for a related problem 
(the Alabama paradox for Hamilton's method).

We use this approach for several different election methods, and calculate,
for given sizes of the parties (under a technical condition, see Sections
\refand{Sresults}{Sdep}), 
the (asymptotic) distribution of the deviation from perfect proportionality
as well as its mean, \ie{} the systematic bias, and its variance.
In particular, this shows how the bias depends on the sizes of the parties in
a more detailed way than previous studies that consider random party sizes 
(see the references given above).

The election methods that we consider are divisor methods of the ``linear''
(or ``stationary'') type (including Jefferson/D'Hondt and Webster/\StL)
and quota methods (including Hamilton/Hare and Droop).
See further  \refS{Snotation}, and Appendices \ref{Adivisor}--\ref{Aquota}.

The main results are stated in \refS{Sresults} and proven in Sections
\ref{SpfD}--\ref{SpfQ}. 
\refS{Srandom} show that the results hold also in the more traditional
approach with deterministic house size $N$ and random party sizes.
The  following sections contain some simple
applications of our main results: 
\refS{Squota} discusses the probability of violating quota,
\refS{Smerge} the expected gain or loss for parties forming 
an \app{} (coalition),
\refS{Sdiv} the \StL{} divergence studied by 
\citet{HPS2004,HPS2005}
and \refS{Sothergood} some related goodness-of-fit functionals.
\refS{Sdep} discusses our technical condition and the case
of rational sizes of the parties.
The appendices contain some background material.

\begin{ack}
  I thank Friedrich Pukelsheim for many helpful comments.
\end{ack}

\section{Notation}\label{Snotation}
We assume throughout that we have $m$ parties with $v_i$ votes for party
$i$; we let $V\=\sumim v_i$ be the total number of votes and $p_i\=v_i/V$
the proportions of votes for party $i$. (In real elections there may 
also be votes for parties that are disqualified because of a threshold rule,
blank votes, and other invalid votes; these are ignored here so $\sumim
p_i=1$. We  consider only parties that have received at least one
vote, so 
$p_i>0$.) We further assume that $N$ seats are to be distributed
(the \emph{house size}), and let 
$s_i$ be the number of seats given to party $i$; thus $\sumim s_i=N$.
We occasionally write $s_i(N)$ when we need to specify the house size.

Strict proportionality would give
\begin{equation}
  \label{q_i}
q_i\=\frac{v_i}VN=p_iN
\end{equation}
seats to party $i$. (This is usually not an integer.) We define the
\emph{seat excess} for party  $i$ to be the difference
\begin{equation}
  \label{gD}
\gD_i\=s_i-q_i=s_i-p_iN.
\end{equation}
Note that
\begin{equation}\label{gD0}
  \sumim \gD_i=\sumim s_i - N = 0.
\end{equation}

We  use the standard notations $\floor x$ and $\ceil x$ for rounding down
and up of a real number $x$, \ie{} the largest
integer $\le x$ and the smallest integer $\ge x$, respectively.
We denote the fractional part of $x$ by $\frax x\=x-\floor x$.

More generally, 
let $\ga$ be a real number. We say that the
\emph{$\ga$-rounding} of a real number $x$ is the integer $\rga x$ such that
\begin{equation}\label{rga}
  x-\ga \le \rga x \le x-\ga+1;
\end{equation}
equivalently,
\begin{equation}\label{rga1}
 \rga x +\ga-1 \le x \le \rga x+\ga.
\end{equation}
If $x-\ga$ is an integer, we regard $\rga x$ as multiple-valued, and accept
both $x-\ga$ and $x-\ga+1$ as values of $\rga x$. (This exceptional
case typically occurs at ties in the election methods.)
Consequently,
\begin{equation}\label{dga1}
  \rga x = \ceil{x-\ga} = \floor{x+1-\ga},
\end{equation}
except the case when $x-\ga$ is an integer and both $x-\ga$ and $x-\ga+1$
are possible values.
If $0\le \ga\le 1$, this means that $x$ is rounded down if its fractional
part is less than $\ga$ and up if its fractional part is greater than $\ga$.
In particular, $\ga=\frac12$ yields standard rounding, $\ga=0$ yields
rounding up and $\ga=1$ yields rounding down.
(But note that we allow also $\ga<0$ or $\ga>1$, in which case $|x-\rga x|$
may be greater than 1.) 

The two types of election methods that we consider can be defined as
follows; see further Appendices \ref{Adivisor}--\ref{Aquota} where
further details and equivalent definitions are given.
Note that all methods that we consider are \emph{homogeneous}, \ie, the
result depends only on the proportions $p_i$ of votes, but not on the total
number $V$ of votes; hence we will mainly use $p_i$ rather than $v_i$ in the
sequel. 

\begin{itemize}
\item 
The \emph{\gbdiv},
or the \emph{divisor method with divisors $d(n)=n+\gb-1$} (where $\gb$ is a
real number): 
Let
\begin{equation}\label{gbdiv}
  s_i\=\rgb{\frac{v_i}D}=\rgb{\frac{p_i}{D'}},
\end{equation}
where $D$ (or $D'=D/V$) is chosen such that $\sumim s_i=N$.
(If $\gb<0$ or $\gb>1$, there are some exceptions for small $N$, see
\refApp{ADlinear}.) 
This includes several important methods, in particular 
$\gb=1$ (Jefferson, D'Hondt),
$\gb=1/2$ (Webster, \StL),
$\gb=0$ (Adams).

\item 
The \emph{\gamq{} method} (where $\gam$ is a real number):
Let $Q\=V/(N+\gam)$ and let
\begin{equation}\label{gamq}
  s_i\=\rga{\frac{v_i}Q}=\rga{(N+\gam)p_i},
\end{equation}
where $\ga$ is chosen such that $\sumim s_i=N$.
In principle, $\gam$ can be any real number (provided $N+\gam>0$); in
practice only the cases $\gam=0$ (\empho{Hamilton},  Hare,
\empho{method  of largest remainder})
and $\gam=1$ (\empho{Droop}) are important,
although also $\gam=2$ has some limited use.
\end{itemize}

\begin{remark}
  \label{Rtie}
In both divisor methods and quota methods ties may occur; then the \gar{} in
\eqref{gbdiv} or \eqref{gamq} has two possible values (because
$x-\ga\in\bbZ$ in \eqref{rga}) for (at least) two parties, with (at least)
one of them rounded to the larger value and (at least) one rounded to the
smaller value. 
We assume that ties are resolved at random, but 
ties will usually  not be a problem for us because of our
assumptions below.
\end{remark}

We let $\U(a,b)$ denote the uniform distribution on the interval $(a,b)$. We
will use $U_i$ and $\tU_i$ for random variables with $U_i\sim\U(0,1)$ and
$\tU_i\sim\U(-\frac12,\frac12)$, and we assume that such variables with
different indices are independent. Recall that then $\E U_i=1/2$,
$\E\tU_i=0$ and $\Var U_i=\Var\tU_i=1/12$.

We denote convergence in distribution by $\dto$,
convergence in probability by $\pto$ 
and
equality in distribution by $\eqd$.

\section{Main results}\label{Sresults}

As explained in the introduction, we assume that $p_1,\dots,p_m$ are fixed,
and let the total number of seats (house size) $N$ be random and uniformly
distributed on \set{1,\dots,\NN} for some large integer $\NN$; we then consider
limits as $\NNtoo$. 
For convenience, we denote this situation by $\Ntoox$.
(Note that this implies $N\pto\infty$.)

\begin{remark}\label{RN}
We assume that $N$ is distributed in this way for simplicity, but the
results can easily be extended to more general  sequence of random
$N\in\bbN$.
For example, we may take $N$ uniformly distributed on an interval
$[N_1(k),N_2(k)]$ and let $\ktoo$, for any two sequences $N_1(k)$ and $N_2(k)$
of positive integers
with $N_2(k)-N_1(k)\to\infty$.
In fact, we conjecture that the results hold for any $N$  
that satisfies the weak condition 
\begin{equation}
  \label{g}
\gf(\theta)
\=\E e^{\ii N\theta}\to0
\qquad\text{for any $\theta\in(0,2\pi)$}.
\end{equation}
(This implies $N\pto\infty$ and excludes parity restrictions.)
Note that Weyl's Lemma \ref{LW} extends to this situation.
However, we have not verified all details for this case, and we leave
further investigations to the interested reader.
\end{remark}

In order to obtain simple asymptotic results, we will usually make one
mathematical simplification:
we assume that $\ppm$ are \emph{\linQ}, \ie, 
that there is no relation
\begin{equation}\label{linq}
a_1p_1+\dots+a_mp_m=0  
\end{equation}
with all coefficients $a_i$ rational and not all $0$.
(Equivalently, there is no relation \eqref{linq} with integer coefficients
$a_i$, not all 0.)
The case when this assumption does not hold is discussed in \refS{Sdep}.

\begin{remark}\label{Rind}
Mathematically, this assumption is reasonable, since if we choose $\ppm$ at
random (uniformly given that their sum is 1, as in \refS{Srandom}), they
will almost surely be \linQ. 
However, for real elections 
the assumption is clearly unreasonable since vote counts $v_i$ are
integers and the $p_i$ thus rational numbers. Nevertheless, 
the results below are good approximations 
if the numbers $p_i$ have large denominators and there are no relations
\eqref{linq} with small integers $a_1,\dots,a_m$,
see \refS{Sdep}. For practical purposes we thus can use the results below as
good approximations for any $p_i$ and any $N$ that is large enough.
(However, we do not investigate the rate of convergence, and how large $N$
has to be.)
\end{remark}

\subsection{Divisor methods}

We begin with a simple deterministic bound, valid for all house sizes and
numbers of votes, \cf{}
\cite[Satz 6.2.11]{Kopfermann}.
It follows from the proof, or
from \refT{TD} below, that the bounds are best possible.
Proofs of the results below are given in \refS{SpfD}.

\begin{theorem}   \label{TDbound} 
Consider the \gbdiv. 
Then
\begin{equation}\label{tdbound1}
p_i-1+\lrpar{\gb-1}(mp_i-1)
\le  \gD_i \le
\lrpar{\gb-1}(mp_i-1) + (m-1)p_i.
\end{equation}
Equivalently,
\begin{equation}\label{tdbound2}
  \lrabs{\gD_i-\Bigpar{\gb-\frac12}(mp_i-1)}
\le \frac12+\frac{m-2}2p_i.
\end{equation}
\end{theorem}

For the asymptotic bias, we have a simple formula.

\begin{theorem}
  \label{TDmean} 
Consider the \gbdiv, and
suppose that 
$p_1,\allowbreak\dots,p_m$ 
are \linQ.
Then, as \Ntoox, 
\begin{equation}
  \label{tdmean}
\E\gD_i \to \lrpar{\gb-\tfrac12}(mp_i-1).	
\end{equation}
\end{theorem}
Note that the asymptotic bias in \eqref{tdmean}
for a party depends only on its size $p_i$
and the number of parties, but not on the distribution of sizes between the
other parties. 
The limit in \eqref{tdmean} agrees with the formula by
\citetq{(S) on p. 208}{Kopfermann} for the advantage (or disadvantage) for
different parties, derived by geometric considerations.

In particular, \refT{TDmean} shows that
the bias is 0 for every party when $\gb=1/2$ (Webster/\StL).
It is well-known that this divisor method is (asymptotically)
unbiased, see \eg{} \cite{Polya1919,BY,Pukel2003b,DSbias2004,DSbias2005} for
various justifications of this; it is satisfying that our approach confirms
this, and shows that the method really is unbiased for a party of any size.

For $\gb>1/2$ (for example Jefferson/D'Hondt with $\gb=1$), the bias is
positive for $p_i>1/m$ and negative for $p_i<1/m$.
It is well-known that these methods favour larger parties; we here see
exactly how much parties of different sizes are favoured, and that the
boundary is exactly at the average size $1/m$.
Note also that the bias grows with $m$, especially for large parties. 
For example, for $\gb=1$, a large
party with $p_i\approx 1/2$ has bias $\approx(m-2)/4$, 
while a small party ($p_i\ll 1/m$) always has a bias $\approx-1/2$.

For $\gb<1/2$ (for example Adams with $\gb=0$), we have the opposite biases;
the boundary is still at $p_i=1/m$.

\begin{example}\label{EStL}
  \citet{StL2} compared in 1910 the results of his method 
(which he called the \emph{method of least squares}) and D'Hondt's method 
in the case of two parties (that he
thought of as two of several parties in an election) with $p_1\ge p_2$.
He did not consider the biases in our sense, comparing the distribution of
the seats to the distribution of the votes, but his approach means studying
the difference of the two biases. (Since we know that \StL's method 
($\gb=1/2$)
is
unbiased, this is equivalent to calculating the bias of D'Hondt's method
($\gb=1$).)

\StL{} found that if $k\=p_2/p_1\in(0,1]$, then the largest party gets on the
average $(1-k)/(2(1+k))$ seats more with D'Hondt's method. This agrees with 
\refT{TDmean}, with $m=2$ and $p_1=1/(1+k)$.

\citet{StL2} also considered the case of random $p_i$, assuming that the
ratio $p_2/p_1$ is uniformly distributed over $(0,1]$, and found by
integration the average gain
\begin{equation}
  \intoi \frac{1-k}{2(1+k)}\dd k = \log 2 - \frac12
\approx 0.193
\end{equation}
for the larger party with D'Hondt's method.
(Note that this differs from the average bias $1/4=0.25$ for the largest
of two parties derived by \cite{Pukel2003b}
assuming that $p_2$ is uniformly distributed on $(0,1/2)$,
see \refT{TPRmean} below with $m=2$.
The reason is that the two probability
distributions of the party sizes are different; they are uniform in
different ways.) 
\end{example}

\begin{remark}\label{RStL}
  \citet{StL2} seems to have thought that this calculation for two parties
  applies to any two of several parties, and that therefore, if party $A$ is
  larger than $B$, then $B$ would on the average gain  $0.193$ seats
  compared to $A$ if the method is changed from D'Hondt's to \StL's,
  regardless of the number and sizes of the other parties. \refT{TDmean}
  shows that this is incorrect; in fact, if $A$ and $B$ both are small
  parties, they both have essentially the same bias $-1/2$ with D'Hondt's
  method. 
It is true that for divisor methods, the distribution of seats to two
parties $A$ and $B$ is the same as if the total number of seats for these
two parties is distributed between them by the same divisor method, see
\refT{Tconsistency}. However, if both $A$ and $B$ are small parties, they
will by \eqref{tdmean} get on the average almost 1 seat less together than
their proportional share, and if their total number is redistributed between
them, most of this loss will be borne by the larger party, which offsets
the advantage of the larger party in the redistribution and explains the
fallacy. 
More generally, parties 1 and 2 will by \eqref{tdmean} on the average get
together $N'\=N(p_1+p_2)+(\gb-1/2)(mp_1+mp_2-2)$ seats. If we distribute this
number of seats between the two parties, the same formula would yield that
party 1 on the average gets, with $p_1'\=p_1/(p_1+p_2)$,
\begin{equation*}
  \begin{split}
  N'p_1'+(\gb-\tfrac12)(2p_1'-1)
&=
Np_1+(\gb-\tfrac12)\bigpar{p_1'(m(p_1+p_2)-2)+2p_1'-1}
\\&
=  Np_1+(\gb-\tfrac12)(mp_1-1)	
  \end{split}
\end{equation*}
seats, in accordance with \refT{TDmean}.
  \end{remark}

We can describe asymptotically not only the mean but also the distribution
of the seat excess. We also obtain a joint limit distribution, which  is
explicit
although a little involved.

\begin{theorem}\label{TD}
Consider the \gbdiv, and
suppose that $\ppm$ are \linQ.
Then, as $\Ntoox$, for each $i$,
\begin{equation}\label{td0}
  \gD_i\dto 
\bX_i\=
\xpar{\gb-\tfrac12}(mp_i-1)
+\tU_0+p_i\sum_{k=1}^{m-2} \tU_k,
\end{equation}
where $\tU_k\sim\U(-\frac12,\frac12)$ are independent.
Moreover, jointly for $i=1,\dots,m$, 
\begin{equation}
\gD_i\dto X_i,   
\end{equation}
where the limit random variables
$X_i$ can be constructed as follows:
Let\/ $U_1,\dots,U_m\sim\U(0,1)$ and let\/ $J\in\set{1,\dots,m}$ be random with
$\P(J=j)=p_j$, $j=1,\dots,m$, with all variables independent. Let
\begin{equation}\label{tdvi}
  V_i\=
\ett{J\neq i}U_i=
  \begin{cases}
	U_i, & i\neq J,\\
0,& i=J,
  \end{cases}
\end{equation}
and let, finally,
\begin{equation}\label{tdxi}
  X_i\=p_i\sumjm V_j-V_i+(\gb-1)(mp_i-1).
\end{equation}
\end{theorem}

Of course, as we shall verify in \refS{SpfD}, $X_i\eqd \bX_i$.

The description in \eqref{tdvi}--\eqref{tdxi} of the joint distribution of
$X_1,\dots,X_m$ is explicit but rather involved, and it is quite possible
that there exists a simpler description; cf.\ \refR{RSA}.
For $m=2$, $X_2=-X_1$, and thus $(X_1,X_2)\eqd(\bX_1,-\bX_1)$ where
$\bX_1\=(2\gb-1)p_1-\gb+U_0$ by \eqref{td0}. We leave it as an open problem
to try to find alternative expressions for $m\ge3$.

\begin{corollary}
  \label{CDvar} 
Consider the \gbdiv, and
suppose that $\ppm$ are \linQ.
Then, as $\Ntoox$, 
  \begin{align}\label{cdvar}
\Var\gD_i &\to \frac{1+(m-2)p_i^2}{12}
\intertext{and, when $i\neq j$,}\label{cdcovar}
\Cov(\gD_i,\gD_j) &\to \frac{(m-2)p_ip_j-p_i-p_j}{12}.
  \end{align}
\end{corollary}

Note that the asymptotic distribution of the seat excess depends on $\gb$
only through its mean $\E X_i$ (\ie, the asymptotic bias in \eqref{tdmean});
the centred random variable $X_i-\E X_i$ does not depend on $\gb$.
We see also, from \eqref{tdxi} or \eqref{cdvar}, that this random part of
the seat excess tends to be larger in absolute value for a large party.

In the limit as $p_i\to0$, \ie{} for very small parties, \eqref{tdvi} and
\eqref{tdxi} yield $V_i\dto U_i$ and 
$$X_i\dto -U_i -(\gb-1)=1-U_i-\gb\eqd U_i-\gb\sim \U(-\gb,1-\gb).$$

\begin{remark}\label{Rgbper}
A simple heuristic motivation for the bias in \refT{TDmean} is that the
$\gb$-rounding in \eqref{gbdiv} is unbiased if $\gb=1/2$, but otherwise, it
costs on the average each party $\gb-1/2$ seats, compared to $v_i/D$ which
is truly proportional. In total, this makes a loss of $m(\gb-1/2)$ seats,
which are redistributed proportionally (by shifting $D$), 
so party $i$ gets back 
$p_im(\gb-1/2)$ seats, on the average, giving a net gain of
$(mp_i-1)(\gb-1/2)$, as shown by \refT{TDmean}.

For a related, but rigorous, argument,
let us compare the methods with parameters $\gb$ and $\gb+1$, temporarily
using a superscript to distinguish the methods. (For example, Adams's method
and Jefferson's method.)
By \refT{Tgb+1} (assuming that $N$ is large enough),
\begin{equation}
  s_i^{(\gb)}(N) = s_i^{(\gb+1)}(N-m)+1.
\end{equation}
Hence,
\begin{equation}\label{rgbper}
  \gD_i^{(\gb)}(N) = s_i^{(\gb)}(N)-Np_i= \gD_i^{(\gb+1)}(N-m)-mp_i+1.
\end{equation}
Taking the limit as \Ntoox,
we see that the asymptotic distributions also differ by $mp_i-1$.
Hence, it is not surprising that the bias in \refT{TDmean} is $mp_i-1$ times a
 linear function of $\gb$.
\end{remark}

\subsection{Quota methods}

We again begin with a simple deterministic bound,
\cf{} \cite[Satz 6.2.3]{Kopfermann}.
It follows from the proof, or
from \refT{TQ} below, that the bounds are best possible.
Proofs of the results below are given in \refS{SpfQ}.

\begin{theorem}\label{TQbound}
For the $\gam$-quota method,
\begin{equation}\label{qbound}
\gam\Bigpar{p_i-\frac1m} - \frac{m-1}m \le
  \gD_i \le \gam\Bigpar{p_i-\frac1m} + \frac{m-1}m.
\end{equation}
\end{theorem}

Again we have a simple formula for the asymptotic bias.
\begin{theorem}\label{TQmean}
Consider the $\gam$-quota method and
suppose that $\ppm$ are \linQ. Then, as \Ntoox,
\begin{equation}\label{tqmean}
\E  \gD_i \to\gam\Bigpar{p_i-\frac1m}.
\end{equation}
\end{theorem}

As for the divisor methods, the asymptotic bias in \eqref{tqmean}
for a party depends only on its size $p_i$
and the number of parties, but not on the distribution of sizes between the
other parties. 
The limit in \eqref{tqmean} agrees with the formula by
\citetq{(S) on p. 199}{Kopfermann} for the (dis)advantage  for
different parties, derived by geometric considerations.

In particular, \refT{TQmean} shows that
the bias is 0 for every party when $\gam=0$ (method of largest
remainder/Hamilton/Hare). 
Again it is well-known that this quota method is (asymptotically)
unbiased, see \eg{} \cite{Polya1919,BY,Pukel2003b,DSbias2004,DSbias2005};
our approach confirms this and shows that the method is unbiased for a
party of any size. 

For $\gam=1$ (Droop), the bias is
positive for $p_i>1/m$ and negative for $p_i<1/m$, 
just as for Jefferson's
method discussed above; hence 
Droop's method too favours larger parties. 
(See \cite[p.\ 118]{Kopfermann} for figures illustrating the cases $m=2$ and
$m=3$.) 
Note, however, that the bias for Droop's method does
not grow with $m$. A comparison between \eqref{tqmean} (with $\gam=1$) and
\eqref{tdmean} (with $\gb=1$) shows that the bias for Jefferson's method is
$m/2$ times the bias for Droop's method, for any party of any size.
Hence  Droop's method is less biased than Jefferson's for every
$m\ge3$. (For $m=2$ not only the biases are the same; in fact the methods
coincide when there are only two parties, see \refApp{SS2}.)

\begin{remark}\label{RQheur}
There is a simple heuristic explanation of the formula \eqref{tqmean} for
the bias, if we accept the fact that Hamilton/Hare's method ($\gam=0$) is
unbiased. Consider $\gam=1$ (Droop); this method can be seen as doing
Hamilton's method with $N+1$ seats, but retracting the last seat. Since
Hamilton's method is unbiased, this first distributes on the average
$(N+1)p_i$ seats to party $i$. The retracted seat is taken essentially
uniformly at random, since it depends on the fractional parts of 
the numbers $q_i=Np_i$ only;
hence each party loses on the average $1/m$ seats in that step. This
argument gives a bias of $p_i-1/m$, in accordance with \eqref{tqmean}.
\end{remark}

Again we have an explicit expression also for the asymptotic distribution
of the seat excess. 

\begin{theorem}\label{TQ}
Consider the $\gam$-quota method and
suppose that $\ppm$ are \linQ. Then, as \Ntoox, for each $i$,
\begin{equation}\label{tq}
  \gD_i\dto 
Y_i\=
\gam\Bigpar{p_i-\frac1m}
+\tU_0+\frac1m\sum_{k=1}^{m-2} \tU_k,
\end{equation}
where $\tU_k\sim\U(-\frac12,\frac12)$ are independent.
\end{theorem}

\begin{corollary}
  \label{CQvar} 
Consider the $\gam$-quota method and
suppose that $\ppm$ are \linQ. Then, as \Ntoox,
\begin{equation}\label{cqvar}
\Var  \gD_i\to 
\frac{1+(m-2)/m^2}{12}
=
\frac{(m+2)(m-1)}{12m^2}.
\end{equation}
\end{corollary}

As for the divisor methods in \refT{TD}, the asymptotic distribution of the
seat excess depends on $\gam$ only through the mean in \eqref{tqmean};
$Y_i-\E Y_i$ does not depend on $\gam$. Moreover, unlike the case of divisor
methods, this centred random variable does not depend on the party size
$p_i$. In particular, for Hamilton/Hare's method ($\gam=0$), the asymptotic
distribution of the seat excess is the same for every party.

\begin{remark}
Let for simplicity $\gam=0$.
  Note that the limit variable $Y_i$ in \eqref{tq}
is uniform $U(-\frac12,\frac12)$ for
  $m=2$; for $m\ge3$ the distribution is more complicated, but as
  $m\to\infty$, the limit converges to $U(-\frac12,\frac12)$ (by \eqref{tq}
and the law of large numbers). Hence, for large $m$, $Y_i$ has essentially
the same distribution as for $m=2$, and although the supremum and infimum of the
range of $Y_i$ tend to 1 and $-1$ as \mtoo, \cf{} \eqref{qbound}, the
probability $\P(Y_i\notin[-\frac12,\frac12])\to0$.

Similarly, it is  seen from the variance formula \eqref{cqvar} that the
variance is 
$1/12\approx 0.08333$ for $m=2$, increases to $5/54\approx0.09259$ for $m=3$
and to the maximum $3/32=0.09375$ for $m=4$; the variance then decreases
again and converges to $1/12$ as \mtoo.

The difference from a uniform variable is not very large, and it may in
practice be a good approximation to regard the seat excess $\gD_i$ in
Hamilton's method as $\U(-\frac12,\frac12)$ for any party and any reasonably
large house size.
\end{remark}

We have stated \refT{TQ} only for individual parties.
It is easy to see that a joint asymptotic distribution exists,
but (unlike the case of divisor methods, \cf{} \refT{TD}) we do not know
any simple description of it.
(There is a description implicit in the proof of \refT{TQjoint},
see \eqref{ain}.)

\begin{theorem}
  \label{TQjoint}
Consider the $\gam$-quota method and
suppose that $\ppm$ are \linQ. Then, as \Ntoox,
$  \gD_i\dto \bY_i$ jointly for $i=1,\dots,m$ for some
random variables $\bY_i$ such that 
the distribution of $(\bY_1-\gam p_1,\dots,\bY_m-\gam p_m)$
is symmetric (exchangeable) and independent of $\ppm$, and
$\bY_i\eqd Y_i$ for each $i$.
\end{theorem}

\begin{corollary}
  \label{CQjoint}
Consider the $\gam$-quota method and
suppose that $\ppm$ are \linQ. Then, as \Ntoox,
when $i\neq j$,
\begin{equation}\label{cqcovar}
  \Cov(\gD_i,\gD_j)
\to
\Cov(\bY_i,\bY_j)=
-\frac{\Var \bY_i}{m-1}
=-\frac{m+2}{12m^2}.
\end{equation}
\end{corollary}

\begin{problem}\label{PQ}
  Find an explicit formula for the limit variables $\bY_1,\dots,\bY_m$ in
  \refT{TQjoint}. 
\end{problem}

\section{Proofs for divisor methods}\label{SpfD}

\begin{proof}[Proof of Theorems \ref{TDbound} and \ref{TD}]
For the \gbdiv,
the seat distribution is given by \eqref{gbdiv} with $D$ and $D'=D/V$
such that $\sumim s_i=N$. We let $D'=1/t$; thus, using \eqref{gbdiv} and
\eqref{dga1}, 
\begin{equation}
  \label{a1}
s_i=\rgb{p_it}=\floor{p_it-\gb+1}
\end{equation}
(with the usual possible exception if $p_it-\gb$ is an integer).
Let here $t$ grow from 0 to $\infty$; 
think of $t$ as time. The seats are assigned one by one (unless there is a
tie).

Fix one party $j$ and let $\tjr$ be the time when party $j$ receives its
$r$:th seat. By \eqref{a1}, $r=p_j\tjr-\gb+1$, and thus 
\begin{equation}
  \tjr=\frac{r+\gb-1}{p_j}.
\end{equation}
At this time $\tjr$, by \eqref{a1} again, party $i$ has $s_i$ seats with
\begin{equation}\label{a1x}
  s_i
=\floor{p_i\tjr-\gb+1}
=\lrfloor{\frac{p_i}{p_j}(r+\gb-1)-\gb+1}
=\lrfloor{\frac{p_i}{p_j}r+\Bigpar{\frac{p_i}{p_j}-1}(\gb-1)}.
\end{equation}
Let, for $i=1,\dots,m$,
\begin{equation}\label{wi}
  W_i
\=\frax{p_i\tjr-\gb+1}
=\lrfrax{\frac{p_i}{p_j}r+\Bigpar{\frac{p_i}{p_j}-1}(\gb-1)}
\in\oi,
\end{equation}
defining $W_i\=1$ in the exceptional case when 
$p_i\tjr-\gb\in\bbZ$ and $s_i=p_i\tjr-\gb$.
(But  $W_i\=0$ when 
$p_i\tjr-\gb\in\bbZ$ and $s_i=p_i\tjr-\gb+1$.)
Then \eqref{a1x} can be written
\begin{equation}\label{a2}
  s_i
=\frac{p_i}{p_j}r+\Bigpar{\frac{p_i}{p_j}-1}(\gb-1)-W_i,
\qquad i=1,\dots,m.
\end{equation}
Note that $W_j=0$, since $s_j=r={p_j\tjr-\gb+1}$.

The total number of seats at the moment $\tjr$ is, by \eqref{a2},
\begin{equation}\label{a2b}
  N=\sumim s_i
=\frac{1}{p_j}r+\Bigpar{\frac{1}{p_j}-m}(\gb-1)-\sumim W_i,
\end{equation}
and thus the seat excess is,
combining \eqref{a2} and \eqref{a2b},
\begin{equation}\label{a3}
\begin{split}
  \gD_i&=  s_i-Np_i
=\Bigpar{\frac{p_i}{p_j}-1}(\gb-1)-W_i
-\Bigpar{\frac{p_i}{p_j}-mp_i}(\gb-1)+p_i\sumkm W_k
\\&
=p_i\sumkm W_k -W_i+(\gb-1)(mp_i-1).
\end{split}
\raisetag{1.5\baselineskip}
\end{equation}
So far everything is deterministic. 
\refT{TDbound} follows from \eqref{a3} since each $W_k\in\oi$ and $W_j=0$.

Let us now assume that $\ppm$ are \linQ; in particular, this implies that
$p_i/p_j$ is irrational when $i\neq j$.
Note that there is a tie between parties $i$ and $j$ at time $t$ if and only
if both $p_it-\gb$ and $p_jt-\gb$ are integers. If parties $i$ and $j$ tie
at two different times $t_1$ and $t_2$, we thus have $p_i(t_1-t_2)\in\bbZ$
and $p_j(t_1-t_2)\in\bbZ$, and taking the quotient we find $p_i/p_j\in\bbQ$,
a contradiction. Hence there is at most one tie for each pair of parties,
and thus at most
$\binom m2$ values of $t$ (and thus of $N$)
for which there is any tie. 
Since we consider asymptotics, we may ignore the ties and assume that no tie
occurs.  

Let now $N$ be random in \set{1,\dots,\NN} and let $J=J(N)$ be the party getting
the $N$:th seat. Then $\gD_i$ is given by \eqref{a3}, where 
$W_i=W_i(N)$ for each $N$ is given by \eqref{wi}  with $j=J$ and $r=s_J(N)$.
Fix again $j$ and consider first only $N$ such that $J=j$. The number of
such $N\le \NN$ is $s_j(\NN)$, the number of seats party $j$ receives when the
house size is $\NN$; note that, \eg{} by \refT{TDbound},
$s_j(\NN)=p_j\NN+O(1)$.  

We recall a  well-known result by Weyl \cite{Weyl}.
(The result is usually stated with $a_1=\dots=a_k=0$; this case implies
immediately the more general version given here.)
We say that an infinite sequence $(v_n)_{n\ge1}\in\oix^{k}$ is
\emph{uniformly distributed} if the empirical distributions
$\nn\qw\sum_{n=1}^{\nn}\gd_{v_n}$
converge to the uniform distribution as \nntoo, where $\gd_{v_n}$ denotes
the Dirac measure.  
(Recall that this means that if $A\subseteq\oix^{k}$ with $\gl(\partial
A)=0$, then 
$\#\set{n\le \nn:v_n\in A}/\nn\to \gl(A)$, 
where $\gl$ is the usual Lebesgue measure, see \eg{} \cite{Billingsley}.)

\begin{lemma}[Weyl]\label{LW}
Suppose that $y_1,\dots,y_k$ and\/ $1$ are \linQ, and let $a_1,\dots,a_k$ be
any real numbers.
Then the sequence of vectors 
$(\frax{ny_1+a_1},  \dots,\frax{ny_k+a_k})_{n\ge 1}\in[0,1)^k$ 
is uniformly distributed in $[0,1)^k$.
\nopf
\end{lemma}
(The standard proof is by showing that the Fourier transform (characteristic
function) 
$\frac1{\nn}\sum_{n=1}^{\nn} 
 \exp\bigpar{2\pi\ii\sum_{j=1}^k \ell_j\frax{ny_j}}\to0$,
as \nntoo, for any fixed integers $\ell_1,\dots,\ell_k$, not all 0,
see for example \cite[Exercises 3.4.2--3]{Grafakos}.)

Since $\ppm$ are \linQ, so are $p_1/p_j,\dots,p_m/p_j$, \ie{}
$\set{p_i/p_j:i\neq j}\cup\set1$. Taking $y_i\=p_i/p_j$ and 
$a_i\=(y_i-1)(\gb-1)$, we thus see from \refL{LW} that $(W_i)_{i\neq j}$,
given by \eqref{wi} for $r=1,\dots,s_j(\NN)$,
are uniformly distributed vectors in $[0,1)^{m-1}$.
In probabilistic notation, regarding $W_i=W_i(N)$ as random variables,
conditioning on $J=j$  and
recalling $W_j=0$ when $J=j$, 
as \NNtoo{} and thus $s_j(\NN)\too$,
\begin{equation}\label{a4}
\bigpar{ (W_1,\dots,W_m)\mid J=j} \dto \bigpar{\ett{i\neq j}U_i}_{i=1}^m.
\end{equation}
Since $\P(J=j)=s_j(\NN)/\NN\to p_j$ as \NNtoo,
the random variable $J$ has asymptotically the distribution given in
\refT{TD}, and it follows from \eqref{a4} that, with $V_i$ defined as in
\refT{TD},  
\begin{equation}
   (W_1,\dots,W_m)\dto(V_1,\dots,V_m).
\end{equation}
Consequently, \eqref{a3} yields $\gD_i\dto X_i$, jointly for all
$i=1,\dots,m$, with $X_i$ defined by \eqref{tdxi}.

To obtain the simpler form \eqref{td0} for an individual $\gD_i$, let
\begin{equation}
  \label{txi}
\tX_i\=X_i-(\gb-\tfrac12)(mp_i-1)=p_i\sumjm V_j-V_i-\tfrac12 mp_i+\tfrac12.
\end{equation}
Then, if $J=i$,
\begin{equation*}
\tX_i 
=
p_i\sum_{j\neq i}U_j -\tfrac12mp_i+\tfrac12  
\eqd p_i\sum_{j=1}^{m-1}\tU_j+\tfrac12-\tfrac12p_i,
\end{equation*}
and if $J\neq i$,
\begin{equation*}
\tX_i 
\eqd
p_i\sum_{j\neq i,J}U_j + (p_i-1)U_i-\tfrac12mp_i+\tfrac12
\eqd p_i\sum_{j=1}^{m-2}\tU_j+(p_i-1)\tU_0-\tfrac12p_i.  
\end{equation*}
Hence, if we define 
\begin{equation}
  \label{zn}
Z\=
\begin{cases}
  p_i\tU_{m-1}+\frac12-\frac12p_i, & J=i, \\
  (p_i-1)\tU_{0}-\frac12p_i, & J\neq i, \\
\end{cases}
\end{equation}
then
\begin{equation*}
  \tX_i\eqd p_i\sum_{j=1}^{m-2}\tU_j+Z,
\end{equation*}
with $\tU_j$ ($1\le j\le m-2$) and $Z$ independent.
By \eqref{zn}, conditioning on $J$, we have the distributions
$(Z\mid J=i)\sim \U(\frac12-p_i,\frac12)$
and
$(Z\mid J\neq i)\sim \U(-\frac12,\frac12-p_i)$; furthermore $\P(J=i)=p_i$,
and it follows that $Z\sim\U(-\frac12,\frac12)$, so $Z\eqd\tU_0$ and thus
\begin{equation*}
  \tX_i\eqd p_i\sum_{j=1}^{m-2}\tU_j+\tU_0.
\end{equation*}
The definition \eqref{txi} now shows that $X_i\eqd\bX_i$, with $\bX_i$
defined in \eqref{td0}.
\end{proof}

\begin{proof}[Proof of  \refT{TDmean} and \refC{CDvar}]
Since $\gD_i$ is bounded, \eg{} by \refT{TDbound} (or by the simpler
\refT{TDprop}), 
the convergence in
distribution $\gD_i\dto X_i$ in \refT{TD} implies convergence of all moments
and mixed moments,
see \cite[Theorem 5.5.9]{Gut}; 
hence it suffices to calculate the expectations,
variances and covariances of $X_i$.

We have $\E U_i=1/2$ and thus by \eqref{tdvi}, 
since $J$ and $U_i$ are independent,
  \begin{equation}\label{evi}
\E V_i=\P(J\neq i)\E U_i=\xpar{1-p_i}\tfrac12.	
  \end{equation}
Consequently, $\E\sumim V_i=\frac12\sumim(1-p_i)=\frac12(m-1)$ and, by
\eqref{tdxi}, 
\begin{equation*}
  \begin{split}
  \E X_i 
&=p_i\E\sumjm V_j-\E V_i+(\gb-1)(mp_i-1)
\\&
= p_i\frac{m-1}2-\frac{1-p_i}2+(\gb-1)(mp_i-1)
=\Bigpar{\gb-\frac12}(mp_i-1).	
  \end{split}
\end{equation*}
This shows \refT{TDmean}.

For the variance and covariances, we have for any $i$
\begin{align*}
\E V_i^2=\P(J\neq i)\E U_i^2=\xpar{1-p_i}\tfrac13	  
\end{align*}
and, when $i\neq j$, by independence,
\begin{align*}
\E (V_iV_j)=\E\bigpar{\ett{J\neq i}\ett{J\neq j}U_iU_j}
=\P(J\notin\set{i,j})\E U_i\E U_j=\frac{1-p_i-p_j}4	  .
\end{align*}
Consequently, using also \eqref{evi},
\begin{align*}
  \Var(V_i)&=\E V_i^2-(\E V_i)^2
=\frac{1-p_i}3-\frac{(1-p_i)^2}4	  
=\frac{1+2p_i-3p_i^2}{12},
\\
\Cov(V_i,V_j)&=
\frac{1-p_i-p_j}4	  
-\frac{1-p_i}2\cdot\frac{1-p_j}2
=-\frac{p_ip_j}4
=-\frac{3p_ip_j}{12}.
\end{align*}
We obtain from this also
\begin{align*}
  \Cov\Bigpar{V_i,\sumjm V_j}
&=
\sumjm  \Cov(V_i, V_j)
=\frac{1+2p_i}{12}-\sumjm\frac{3p_ip_j}{12}
\\&
=\frac{1+2p_i}{12}-\frac{3p_i}{12}
=\frac{1-p_i}{12},
\\
\Var\Bigpar{\sumim V_i}
&= \sumim \Cov\Bigpar{V_i,\sumjm V_j}
= \sumim \frac{1-p_i}{12}
=\frac{m-1}{12}.
\end{align*}
Consequently,
\eqref{tdxi} yields
\begin{align*}
  \Var(X_i)
&=
p_i^2\Var\Bigpar{\sumjm V_j} -2p_i  \Cov\Bigpar{V_i,\sumjm V_j}
+\Var(V_i)
\\&
=p_i^2\frac{m-1}{12}-2p_i\frac{1-p_i}{12}+\frac{1+2p_i-3p_i^2}{12}
=\frac{(m-2)p_i^2+1}{12}
\intertext{and, for $i\neq j$,}
\Cov(X_i,X_j)&=
p_ip_j\Var\Bigpar{\sumkm V_k} -p_i  \Cov\Bigpar{V_j,\sumkm V_k}
-p_j  \Cov\Bigpar{V_i,\sumkm V_k}
\\&\hskip6em
+\Cov(V_i,V_j)
\\&
=p_ip_j\frac{m-1}{12}-\frac{p_i(1-p_j)+p_j(1-p_i)}{12}-\frac{3p_ip_j}{12}
\\&
=\frac{(m-2)p_ip_j-p_i-p_j}{12}.
\end{align*}
\refC{CDvar} follows.
\end{proof}

\section{Proofs for quota methods}\label{SpfQ}

\begin{proof}[Proof of Theorems \refand{TQbound}{TQ}]
The proof for quota methods uses arguments similar to the proof in
\cite{SJ258}, and the reader might wish to compare the versions. 
(Only the case $\gam=0$ is treated in \cite{SJ258}, but that is only a minor
simplification.
For completeness we repeat some results from \cite{SJ258}.)
For notational convenience, we show the results for party 1, \ie{} we take
$i=1$ in the proofs. 

By definition, $s_j=\rga{(N+\gam)p_j}$ for an $\ga$ that makes $\sumjm s_j=N$.
Explicitly, by \eqref{rga1},
\begin{equation}\label{q1}
  s_j+\ga-1\le(N+\gam)p_j\le s_j+\ga.
\end{equation}
Thus, using party 1 as a benchmark,
\begin{equation}\label{q2}
  s_j-s_1-1\le(N+\gam)(p_j-p_1)\le s_j-s_1+1.
\end{equation}

Since $\gD_j-\gD_1=s_j-s_1-N(p_j-p_1)$, \eqref{q2} yields
\begin{equation}
  \gD_j-\gD_1-1\le\gam(p_j-p_1)\le \gD_j-\gD_1+1.
\end{equation}
Thus, summing over all $j\neq 1$, recalling that $\sumjm\gD_j=0$ and $\sumjm
p_j=1$, 
\begin{equation}
  -m\gD_1-(m-1)\le\gam(1-mp_1)\le -m\gD_1+(m-1).
\end{equation}
Hence, 
\begin{equation}
\gam\Bigpar{p_1-\frac1m} - \frac{m-1}m \le
  \gD_1 \le \gam\Bigpar{p_1-\frac1m} + \frac{m-1}m, 
\end{equation}
which shows \refT{TQbound}. 

Suppose now that all $p_i-p_j$ are irrational. A tie between parties $i$ and
$j$ can occur only if $(N+\gam)p_i-(N+\gam)p_j\in\bbZ$.
If this happens for two different house sizes $N_1$ and $N_2$, then both
$(N_1+\gam)(p_i-p_j)\in\bbZ$ and
$(N_2+\gam)(p_i-p_j)\in\bbZ$, and thus by taking the difference
$(N_1-N_2)(p_i-p_j)\in\bbZ$, which is impossible if
$p_i-p_j\notin\bbQ$. Hence, our assumption implies that there is a tie
between parties $i$ and $j$ for at most one value of $N$, and thus at most
$\binom m2$ values of $N$ for which there is any tie. Since we consider
asymptotics, we can ignore these $N$ and assume that there are no ties.
Hence we can choose $\ga$ such that there are strict inequalities in
\eqref{q1}, and thus there are strict inequalities in \eqref{q2}.

By \eqref{q2} (with strict inequalities)
$\floor{(N+\gam)(p_j-p_1)}\in\set{s_j-s_1-1,\,s_j-s_1}$.
Let 
\begin{equation}\label{ii}
I_j\=
s_j-s_1-\floor{(N+\gam)(p_j-p_1)}.   
\end{equation}
Then $I_j\in\set{0,1}$ and $I_1=0$.
We have
\begin{equation}
  \begin{split}
\gD_j-\gD_1
&
=s_j-s_1-N(p_j-p_1)
\\&
=s_j-s_1-(N+\gam)(p_j-p_1)+\gam(p_j-p_1)
\\&
=I_j-\frax{(N+\gam)(p_j-p_1)}+\gam(p_j-p_1)	.
  \end{split}
\end{equation}
Summing over $j$ we obtain, again recalling $\sumjm \gD_j=0$,
\begin{equation}\label{sw}
  \begin{split}
-m\gD_1=	
\sumjm I_j-\sumjm \frax{(N+\gam)(p_j-p_1)}+\gam(1-mp_1)	.
  \end{split}
\end{equation}
Let $L\=\sumjm I_j$. 
By \eqref{sw}, we then have the formula
\begin{equation}\label{qk}
\gD_1=	
\gam\Bigpar{p_1-\frac1m}
+\frac1m\lrpar{\sumjjm \frax{(N+\gam)(p_j-p_1)}-L}.
\end{equation}
Note that $L$ is an integer with  $0\le L\le m-1$ and, by  \eqref{ii},
\begin{equation*}
  \begin{split}
L&\=\sumjm I_j
  =\sumjm s_j- m s_1-\sumjm \floor{(N+\gam)(p_j-p_1)}
\\&\phantom:
  =N- m s_1-\sumjjm \floor{(N+\gam)(p_j-p_1)}
\\&\phantom:
  \equiv N-\sumjjm \floor{(N+\gam)(p_j-p_1)} \pmod m.	
  \end{split}
\end{equation*}
Let $\Modm(x)\=m\frax{x/m}$, the remainder when $x$ is
divided by $m$. We thus have 
\begin{equation}\label{ql}
  L=
  \Modm\biggpar{ N-\sumjjm \floor{(N+\gam)(p_j-p_1)}}.
\end{equation}

Let $y_j\=p_j-p_1$ and $a_j\=\gam(p_j-p_1)$; thus $(N+\gam)(p_j-p_1)=Ny_j+a_j$.
Let us now assume that \ppm{} are \linQ.
It is easily seen that then 
$y_2,\dots,y_m$ and $\sum_1^m p_j=1$ also are \linQ.
(This can be seen as a change of basis, using a non-singular integer matrix,
in a vector space of dimension $m$ over $\bbQ$.)

We need the following extension of \refL{LW}, slightly generalizing 
\cite[Lemma  4.3]{SJ258}. The proof is given in \refApp{Auni}.

\begin{lemma}\label{Lux}
Suppose that $y_1,\dots,y_k$ and\/ $1$ are \linQ, and let
$a_1,\dots,a_k$ be any real numbers.

Let\/ $\ell_n=
\Modm\bigpar{n- \sumjk\floor{ny_j+a_j}}\in\set{0,\dots,m-1}$. 
Then the sequence of vectors 
\begin{equation}
\label{qulu}
(\frax{ny_1+a_1},  \dots,\frax{ny_k+a_k}, \ell_n)
\in[0,1)^k\times\set{0,\dots,m-1}
\end{equation}
is uniformly distributed in $[0,1)^k\times\set{0,\dots,m-1}$.
\end{lemma}

Note that $L$ is \eqref{ql} equals
$\ell_N$ in \refL{Lux}, with $y_j$ and $a_j$ as chosen above. 
Consequently, \eqref{qk} and \refL{Lux} show that
\begin{equation}\label{qu}
  \gD_1\dto 
\gam\Bigpar{p_1-\frac1m}
+\frac1m\lrpar{\sumjjm U_j-\Lx},
\end{equation}
where $U_j\sim\U(0,1)$ and $\Lx$ is uniformly distributed on
\set{0,\dots,m-1}, with all variable independent.
Moreover, $U_2\eqd1-U_2$ and $\Lx+U_2\sim\U(0,m)$, and thus
\begin{equation*}
\Lx-U_2\eqd \Lx+U_2-1\eqd mU_0-1
\eqd m\tU_0+(m-2)/2.
\end{equation*}
Hence, \eqref{qu} can be written
\begin{equation*}
  \gD_1\dto 
\gam\Bigpar{p_1-\frac1m}
+\frac1m\lrpar{\sum_{j=3}^{m} \tU_j-m\tU_0},
\end{equation*}
with $\tU_j\=U_j-\frac12$ for $j\ge1$,
which by $\tU_0\eqd-\tU_0$ proves \refT{TQ}.
\end{proof}

\begin{proof}[Proof of  \refT{TQmean} and \refC{CQvar}]
Since $\gD_i$ is bounded, \eg{} by \refT{TQbound},
the convergence in
distribution $\gD_i\dto Y_i$ in \refT{TQ} implies convergence of all moments
and it suffices to calculate the expectation and variance of $Y_i$.
This is immediate since $\E\tU_k=0$ and $\Var\tU_k=1/12$; thus
by \eqref{tq}, 
\begin{align*}
  \E Y_i &=\gam\Bigpar{p_i-\frac1m},
\\
\Var Y_i &= \frac1{12}+\frac{m-2}{12 m^2}.
\end{align*}
\refT{TQmean} and \refC{CQvar} follow.
\end{proof}

\begin{proof}[Proof of \refT{TQjoint}]
  Let $W_i\=\frax{(N+\gam)p_i}$.
For any $j$, $\set{p_i}_{i\neq j}$ and $1=\sumim p_i$ are \linQ, and thus
the vectors $(W_i)_{i\neq j}$ are uniformly distributed in $\oix^{m-1}$ by
\refL{LW}. 
Moreover, 
\begin{equation*}
\sumim W_i\equiv\sumim(N+\gam)p_i=N+\gam\equiv \gam \pmod1,  
\end{equation*}
so  $W_j\equiv \gam-\sum_{i\neq j} W_i\pmod 1$ and
\begin{equation*}
  W_j=\Bigfrax{\gam-\sum_{i\neq j} W_i}.
\end{equation*}
Consequently, $(W_1,\dots,W_m)\dto(V_1,\dots,V_m)$ where each
$V_i\sim\U(0,1)$, any $m-1$ of them are independent, and any such set
determines the last variable by $V_j=\frax{\gam-\sum_{i\neq j} V_i}$.
Obviously, the distribution of $(V_1,\dots,V_m)$ is symmetric.

We have
\begin{equation}\label{ai}
  \gD_i-\gam p_i
= s_i-(N+\gam)p_i
= s_i-\floor{(N+\gam)p_i}-W_i.
\end{equation}
By the construction of the \gamqm, the numbers $s_i-\floor{(N+\gam)p_i}$
depend only on $(W_1,\dots,W_m)$, say 
$s_i-\floor{(N+\gam)p_i}=f_i(W_1,\dots,W_m)$ for some function $f_i$, and it
follows from \eqref{ai} that, jointly for all $i$,
\begin{equation}\label{ain}
  \gD_i-\gam p_i \dto Z_i\=f_i(V_1,\dots,V_m)-V_i,
\end{equation}
where the distribution of $(Z_1,\dots,Z_m)$ is symmetric. 
(The $m$ functions $f_i$ differ only by permutations of the variables.)
Now let
$\bY_i\=Z_i+\gam p_i$.
Then $\gD_i\dto\bY_i$, and thus $\bY_i\eqd Y_i$ by comparison with \refT{TQ}.
\end{proof}

\begin{proof}[Proof of \refC{CQjoint}]
By \refT{TQjoint}
and the fact that each $\gD_i$ is boun\-ded, 
$\Cov(\gD_i,\gD_j)\to\Cov(\bY_i,\bY_j)$.
Moreover, by symmetry, $\Cov(\bY_i,\bY_j)$ is the same for all pairs $(i,j)$
with $i\neq j$. 
We have $\sumim \bY_i=0$ as a consequence of
\eqref{gD0}, and thus
\begin{equation*}
  0=\Cov\Bigpar{\bY_i,\sumjm \bY_j}
=\Var(\bY_i)+\sum_{j\neq i}\Cov(\bY_i,\bY_j);
\end{equation*}
hence by symmetry $\Cov(\bY_i,\bY_j)=-\Var(\bY_i)/(m-1)=-\Var(Y_i)/(m-1)$ for
$i\neq j$. 
The final equality follows from \eqref{cqvar}.
\end{proof}

\section{Random party sizes}\label{Srandom}

As said in the introduction, there is a long tradition of studying the bias
of election methods theoretically with a fixed house size and random party
sizes.
More precisely, in this approach one (usually)
assumes that the vector $(\ppm)$ is
random and uniformly distributed over the simplex
$\sxm\=\set{(\ppm)\in(0,1)^m:\sumim p_i=1}$ for some given $m$.
One then orders the party sizes $\ppm$ as $p\oy 1\ge\dots\ge p\oy m$ and
considers 
the seat excesses $\gD\oy1,\dots,\gD\oy m$ of the parties in the same order.
(Thus, $\gD\oy1$ is the seat excess of the largest party, and so on.)
Results for the mean $\E\gD\oy k$ of the seat excess for the $k$:th largest
party (either exact formulas for a fixed $N$ or asymtotic formulas as \Ntoo)
have been obtained by 
\citet{Polya1918,Polya1918b,Polya1919},
\citet{Pukel2003b},
\citet{HPS2005}
and 
Drton and Schwingenschlögl \cite{DSbias2004,DSbias2005};
results on the variance are given by \citet{DSvariance}.
(We follow these papers and order the parties in
decreasing order, \cf{} \refApp{Apox}.)

\citet{HPS2004,HPS2005} and \citet{Sch2008}
have also, more generally, considered other
distributions of $\ppm$ (with an absolutely continuous distribution on $\sxm$).
Note also that \citet{StL2} considered a different distribution for two
parties, see \refE{EStL} (with $p_2/p_1$ uniformly distributed, which is not
the same as $p_2$ uniform).

We cannot obtain these results for a given $N$ directly
from our results for random $N$,
nor can we obtain our results for fixed $\ppm$ from these results.
However, by combining the methods, we can translate one type of results to the
other and thus see that our main results transfer to the case of random $\ppm$.
(Now with random $p_i$ in \eg{} \eqref{td0} and \eqref{tq}.)

\begin{theorem}\label{Trand}
Let $\ppm$ be random with an absolutely continuous distribution on 
the simplex $\sxm$,
and let \Ntoo{} (with $N$ deterministic). 
  \begin{romenumerate}
\item 
For the \gbdiv, the conclusions \eqref{td0}--\eqref{tdxi} of \refT{TD}
hold, with $\tU_i$ and $U_i$ independent of $\ppm$ and with
$\P(J=j\mid\ppm)=p_j$. 

\item 
For the \gamqm, the conclusions of Theorems \refand{TQ}{TQjoint} hold; in
particular \eqref{tq} holds, 
with $\tU_k\sim\U(-\frac12,\frac12)$  independent of each other and of \ppm.
\end{romenumerate}
\end{theorem}

\begin{proof}
Consider first the \gbdiv.
  The proof of \cite[Theorem 3.1(i)]{HPS2005} shows that as \Ntoo,
  \begin{equation}
	\label{qa}
(\gD_1,\dots,\gD_m)\dto(Z_1,\dots,Z_m)
  \end{equation}
for some random variables $Z_j$; with the notation in \cite{HPS2005},
\begin{equation}
  Z_j=-\bigpar{U_j+q-\tfrac12+\sgn(D_c)m_j(D_c)}+c(q-\tfrac12)W_j.
\end{equation}
Here $W_j$ is our $p_j$, $c$ is our $m$, $q$ is our $\gb$, $U_j$ is our
$\tU_j$ and $\sgn(D_c)m_j(D_c)$ is an explicit but rather complicated
function of $U_1,\dots,U_m$ and $W_1,\dots,W_m$;
however, the point is that for
us it suffices to know the existence of some variables 
$Z_j$ such that \eqref{qa}
holds; we can ignore the construction.
The proof of \eqref{qa}
is based on the fact that for any (deterministic) sequence
$\nu_n\to\infty$, the vectors of fractional parts 
$(\frax{\nu_n p_1},\dots,\frax{\nu_n p_{m-1}})$ are uniformly distributed in 
$\oix^{m-1}$, see \refL{Luni}, which is the analogue of \refL{LW} in the
present context.

Consider now $N$ uniformly random on \set{1,\dots,\NN} and let \NNtoo{},
\ie, suppose $\Ntoox$ in the notation used in the previous sections.
Then $N\pto\infty$, so by conditioning on $N$, it follows
that \eqref{qa} holds for such random $N$ as well.
On the other hand, since the distribution of $\ppm$ is absolutely
continuous, almost surely $\ppm$ are \linQ. Thus, by conditioning on $\ppm$
we may apply \refT{TD}, which shows that 
  \begin{equation}
	\label{kni}
(\gD_1,\dots,\gD_m)\dto(X_1,\dots,X_m),
  \end{equation}
where $X_1,\dots,X_m$ are given by \eqref{tdvi}--\eqref{tdxi} (with our
random \ppm).
In the case of random \ppm{} and random $\Ntoox$, we thus have both \eqref{qa}
and \eqref{kni}; consequently,
\begin{equation}
  \label{mar}
(Z_1,\dots,Z_m)\eqd(X_1,\dots,X_m).
\end{equation}

Consider again deterministic \Ntoo. We know that \eqref{qa} holds, so
\eqref{mar} implies that \eqref{kni} holds in this case too. Finally,
$(\bX_i\mid p_i)\eqd(X_i\mid p_i)$ by \refT{TD}, and thus $\bX_i\eqd X_i$, so
\eqref{td0} too holds in the present setting.

For the \gamqm, the proof of \refT{TQjoint} applies to the present
situation as well, provided we use \refL{Luni} instead of \refL{LW}.
In particular, \eqref{ain} holds, with
$V_1,\dots,V_m\sim\U(0,1)$ independent of \ppm. 
Further. 
$(\bY_i\mid\ppm)=(Y_i\mid\ppm)$, and thus $\bY_i\eqd Y_i$ in the case of
random \ppm{} too.  The results follow.
\end{proof}

\begin{remark}
  By considering $\ppm$ uniformly distributed on a small neighbourhood
  $\max_i|p_i-p^0_i|<\eps$ of some $(p_1^0,\dots,p_m^0)\in\sxm$, and then
  letting $\eps\to0$, we recover the limit distributions in \refS{Sresults},
  with a somewhat different interpretation as a double limit. (In this
  version $\ppm$ can be any numbers, including rationals.)
\end{remark}

In particular, we can apply \refT{Trand} to the case of uniformly
distributed \ppm{} taken in decreasing order as above, since then
$(p\oy1,\dots,p\oy m)$ has a uniform distribution on the subset 
$\sxmge\=\set{(\ppm)\in\sxm:p_1\ge\dots\ge p_m}$ of $\sxm$.
This gives immediately the following theorem.
(We leave the case of joint distributions  to the reader.)

\begin{theorem}\label{TPR}
  Let $\ppm$ be random and uniformly distributed on $\sxm$, and let \Ntoo.
  \begin{romenumerate}
\item 
For the \gbdiv, for each $j\le m$,
\begin{equation}\label{trpd}
 \gD\oy j \dto \hX_j 
\=  
\xpar{\gb-\tfrac12}(mp\oy j-1)
+\tU_0+p\oy j\sum_{k=1}^{m-2} \tU_k.
\end{equation}

\item 
For the \gamqm, for each $j\le m$,
\begin{equation}\label{trpq}
  \gD\oy j\dto 
\hY_j\=
\gam\Bigpar{p\oy j-\frac1m}
+\tU_0+\frac1m\sum_{k=1}^{m-2} \tU_k.
\end{equation}
  \end{romenumerate}
Here $\tU_k\sim\U(-\frac12,\frac12)$ are independent of each other and of \ppm.
\nopf
\end{theorem}

As before, this theorem implies moment convergence because $\gD\oy j$ is
bounded. In particular, for the bias in this approach we obtain the
following result. This was (except for the case $\gam\neq0$ in (ii))
conjectured by \citet{Pukel2003b} and proven by \citet{DSbias2005} and
\citet{HPS2005}. 

\begin{theorem}\label{TPRmean}
  Let $\ppm$ be random and uniformly distributed on $\sxm$, and let \Ntoo.
  \begin{romenumerate}
\item 
For the \gbdiv, for each $j\le m$,
\begin{equation}\label{trpdmean}
\E \gD\oy j \to \E\hX_j 
=  
\xpar{\gb-\tfrac12}(m\E p\oy j-1)
=\xpar{\gb-\tfrac12}\biggpar{\sum_{i=j}^m\frac1i-1}.
\end{equation}

\item 
For the \gamqm, for each $j\le m$,
\begin{equation}\label{trpqmean}
\E  \gD\oy j\to 
\E\hY_j=
\gam\Bigpar{\E p\oy j-\frac1m}
=\frac{\gam}m \biggpar{\sum_{i=j}^m\frac1i-1}.
\end{equation}
  \end{romenumerate}
\end{theorem}

\begin{proof}
  The convergence to $\E\hX_j$ or $\E\hY_j$ follows, as just said, from
  \refT{TPR}. The first equalities in \eqref{trpdmean} and \eqref{trpqmean}
  follow immediately from \eqref{trpd} and \eqref{trpq} since $\E\tU_k=0$,
  \cf{} \eqref{tdmean} and \eqref{tqmean}.
It remains only to compute $\E p\oy j$. This is a standard formula, and for
completeness we give a proof in \refApp{Apox},
see \eqref{epox}.
\end{proof}

We similarly obtain results for the variance.
This was found for $m=2$ and $3$ (except the case $\gam\neq0$) by \citet{DSvariance}.

\begin{theorem}\label{TPRvar}
  Let $\ppm$ be random and uniformly distributed on $\sxm$, and let \Ntoo.
  \begin{romenumerate}
\item 
For the \gbdiv, for each $j\le m$,
\begin{equation}\label{trpdvar}
  \begin{split}
&\Var \gD\oy j \to \Var \hX_j 
= \frac1{12}\bigpar{1+(m-2)\E p\oy j^2}
  +(\gb-\tfrac12)^2m^2\Var (p\oy j)	
\\&\hskip5em
=\frac1{12}+
\frac{m-2}{12m(m+1)}\sum_{i=j}^m\frac1{i^2} 
+\frac{m-2}{12m(m+1)}\biggpar{\sum_{i=j}^m\frac1i}^2
\\&\hskip 7em
+(\gb-\tfrac12)^2\lrpar{\frac m{m+1}\sum_{i=j}^m\frac1{i^2} 
-\frac1{m+1}\biggpar{\sum_{i=j}^m\frac1i}^2
}
.
  \end{split}
\end{equation}

\item 
For the \gamqm, for each $j\le m$,
\begin{equation}\label{trpqvar}
  \begin{split}
&\Var  \gD\oy j\to 
\Var\hY_j
= \frac1{12}\Bigpar{1+\frac{m-2}{m^2}}
  +\gam^2\Var (p\oy j)
\\
&=
 \frac{(m+2)(m-1)}{12m^2}
+\gam^2\lrpar{\frac1{m(m+1)}\sum_{i=j}^m\frac1{i^2} 
-\frac1{m^2(m+1)}\biggpar{\sum_{i=j}^m\frac1i}^2
}
.	
  \end{split}
\end{equation}
  \end{romenumerate}
\end{theorem}

\begin{proof}
  Again, by \refT{TPR}, we only have to calculate $\Var \hX_j$ and
  $\Var\hY_j$; the calculations are similar and
we give the details only for $\hX_j$.
We use the standard decomposition
\begin{equation}\label{anal}
  \Var \hX_j
=
\E\bigpar{\Var(\hX_j\mid\ppm)}
+
\Var\bigpar{\E(\hX_j\mid\ppm)}.
\end{equation}
The conditional variance and mean are obtained immediately from
\eqref{trpd}, which yields, \cf{} \eqref{tdmean} and \eqref{cdvar},
\begin{align}
  \E(\hX_j\mid\ppm) &= \xpar{\gb-\tfrac12}(mp\oy j-1),
\\
  \Var(\hX_j\mid\ppm) &= \frac{1+(m-2)p\oy j^2}{12}. \label{f2g}
\end{align}
The first equality in \eqref{trpdvar} follows from \eqref{anal}--\eqref{f2g}.
It remains only to find $\E p\oy j^2$ and $\Var (p\oy j)$; these are
well-known and we give a computation in \refApp{Apox} for completeness, see
\eqref{epox}--\eqref{varpox}. 
\end{proof}

\begin{example}
  For $m=3$ and Jefferson/D'Hondt's method we obtain from \eqref{trpdmean} with
  $\gb=1$,
$\E(\gD\oy 1,\gD\oy 2,\gD\oy 3)\to (\frac5{12},-\frac1{12},-\frac4{12})$,
as found already by \citet{Polya1919}.

For $m=3$ and Droop's method we obtain from \eqref{trpqmean} with
  $\gam=1$,
$\E(\gD\oy 1,\gD\oy 2,\gD\oy 3)\to (\frac5{18},-\frac1{18},-\frac4{18})$.
\end{example}

\begin{example}
  For $m=3$, we have for the \gbdiv
  \begin{multline*}
(\Var\gD\oy1, \Var\gD\oy2, \Var\gD\oy3)
\to\\
\lrpar{
\frac{301}{2592}+\frac{13}{72}(\gb-\tfrac12)^2,
\frac{235}{2592}+\frac{7}{72}(\gb-\tfrac12)^2,
\frac{220}{2592}+\frac{4}{72}(\gb-\tfrac12)^2
}
  \end{multline*}
and for the \gamqm
  \begin{multline*}
(\Var\gD\oy1, \Var\gD\oy2, \Var\gD\oy3)
\to
\\
\lrpar{
\frac{5}{54}+\frac{13}{648}\gam^2,
\frac{5}{54}+\frac{7}{648}\gam^2,
\frac{5}{54}+\frac{4}{648}\gam^2
}.
  \end{multline*}
These were found by
\citet{DSvariance} (where, however, the case $\gam\neq0$ is not treated).
\end{example}

Asymptotic covariances can be computed in the same way, using
\eqref{cdcovar}, \eqref{cqcovar} and \eqref{covpox}.

  In the same way as in Theorems \ref{TPRmean}--\ref{TPRvar}, we can obtain
  results for the mean and variance of $\gD\oy j$ for
  random $\ppm$ with other distributions than the uniform one.
In particular, this gives the formulas by \citet{HPS2005} 
and \citet{Sch2008} for the asymptotic bias in terms of $\E p\oy
j=\E(p_j\mid p_1\ge \dots\ge p_m)$ 
(assuming that $(\ppm)$ has a symmetric distribution). 

For example,
\citet{HPS2005} consider, as a special case, 
the case of a threshold $t$ (with $0<t<1/m$)
and uniform
distribution on $p_1\ge\dots\ge p_m\ge t$. In this case, one has (as in
\cite{HPS2005}) 
$(p_1-1/m,\dots,p_m-1/m)\eqd (1-mt)(p^*_1-1/m,\dots,p^*_m-1/m)$,
where $(p^*_1,\dots,p^*_m)$ has a uniform distribution on $\sxmge$, and
the mean and variance are easily computed.

\section{The probability of violating quota}\label{Squota}

We say that a seat assignment $s_i$ 
\emph{satisfies lower quota} if $s_i\ge\floor {q_i}$ 
and \emph{satisfies upper quota} if $s_i\le\ceil{q_i}$; it \emph{satisfies
  quota} if both hold. In terms of the seat excess $\gD_i=s_i-q_i$,
the assignment satisfies lower [upper] quota if and only if $\gD_i>-1$
[$\gD_i<1$], and it satisfies quota if and only if $|\gD_i|<1$. 

As is well-known, Hamilton/Hare's method always satisfies quota,
while Jefferson's and Droop's methods satisfy lower quota and Adams method
satisfies 
upper quota, see \eg{} \cite{BY} or Theorems \refand{TDbound}{TQbound}
above. 
It is also well-known that Webster/\StL{} does not always satisfy quota, but
that violations are unusual in practice.

Theorems \refand{TD}{TQ} enable us to calculate the (asymptotic)
probabilities that quota is violated for various methods.
We give two examples.

\begin{example} 
  Consider Jefferson/D'Hondt's method.
By \refT{TDbound} with $\gb=1$, $p_i-1\le\gD_i\le(m-1)p_i$; hence lower
quota is always satisfied (as is well-known \cite{BY}),
but upper quota can be violated for a party with $p_i\ge1/(m-1)$.
Consider the case $m=3$, and suppose that $p_1\ge1/2$.
(Then $p_2,p_3<1/2$, so quota is always satisfied for the other two parties;
thus at most one party can violate quota.)
By \refT{TD} with $\gb=1$, letting $\tU_i=\frac12-U_i$,
\begin{equation}
  \gD_1\dto\bX_1=\tfrac12(3p_1-1)+\tU_0+p_1\tU_1
=2p_1-U_0-p_1U_1;
\end{equation}
hence the asymptotic probability of violating (upper) quota is
\begin{equation}
  \P(\bX_1\ge1)=\P\bigpar{U_0+p_1U_1<2p_1-1}=\frac{(2p_1-1)^2}{2p_1}.
\end{equation}
(The set of allowed $(U_0,U_1)$ is a right-angled triangle with sides $2p_1-1$
and $(2p_1-1)/p_1$.)

If we let $(p_1,p_2,p_3)$ be random and uniformly distributed on $\sxx3$, as
in \refS{Srandom}, then $p_1$ has the density $2(1-p_1)$ and thus the
asymptotic probability that party 1 violates (upper) quota is
\begin{equation}
  \int_{1/2}^1 \frac{(2p-1)^2}{2p} \cdot 2(1-p) \dd p
=\ln 2 -2/3 \approx 0.026.
\end{equation}
Consequently, the asymptotic probability that some party violates (upper)
quota is 
\begin{equation}
 3\ln 2 -2 \approx 0.079,
\end{equation}
as found by \citet{NieW}.
\end{example}

\begin{example}
  Consider Jefferson/D'Hondt's method, and a party $i$ with three times the
  average size: $p_i=3/m$. Then the bias $\E\bX_i=1$, by \refT{TDmean} with
  $\gb=1$. It follows by \eqref{td0} and symmetry that $\P(\bX_i>1)=1/2$, so
  the (asymptotic) probability that the party violates quota is $1/2$.
For a larger party, the probability is even greater.
\end{example}

\begin{example}
The Swedish parliament contains 
after the general election in 2010 8 parties: two large with 30\% of the
votes each and 6 small with 5--8\% percent each. 
The seats are in principle distributed by \StL's method. (We ignore here
complications due to the division into 29 constituencies and the system with
adjustment seats, which in 2010 did not give complete adjustment because
the number of adjustment seats was insufficient.)

We note from \refT{TDbound} (with $\gb=\frac12$) that the small parties
always satisfy quota. In fact, \eqref{tdbound2} shows that for Webster/\StL,
only parties with $p_i\ge 1/(m-2)$ can violate quota.

For the large parties we have (approximatively) $p_i= 0.3$, and thus by
\eqref{td0} 
(still with $\gb=\frac12$)
$\bX_i=\tU_0+0.3\sum_{k=1}^6\tU_k$. 
An integration (using Maple) yields
$\P(\bX_i\ge1)=\P(\bX_i\le-1)\approx0.00045$. 
Hence, for each of the two large parties 
 the (asymptotic) probability of violating quota 
is  $0.0009$.
\end{example}

\section{Apparentements}\label{Smerge}  

Suppose that two or more parties decide to form a coalition 
in the election, so that their votes are counted together.
(In some election systems, parties can register such a coalition, called an
\app, and continue to have separate names and lists. Otherwise, the parties
can always appear on election day as one party with a common umbrella name.
In any case we assume that this does not attract or repel any voters, so the
coalition 
gets exactly as many votes as the parties would have had separately, and
that all other parties remain unaffected with the same number of votes as
before.) 

It is well-known that for a \gbdiv{}  with $\gb\ge1$ 
(\eg{} Jefferson/D'Hondt's method) parties can never lose by
forming an \app{}; they will get at least as many seats together
as if they appear
separately. When $\gb\le0$ (\eg{} Adams's method), the opposite is
true, and parties will never 
gain by forming an \app{} (conversely, they may gain by splintering),
and for $0<\gb<1$ (\eg{} Webster/\StL's method), they may both gain and lose.
For D'Hondt's method, this effect is politically important in reality,
especially in small constituencies, and \app{s} are or have been a regular
feature in 
many countries, see \eg{} \cite{Carstairs}.

\citet{Pukel2009b} give a detailed study of the resulting gains in a set of
real elections, including examples, statistics and theoretical values 
assuming that the party sizes are random (as in \refS{Srandom}).
We can now give similar theoretical results for given party sizes.
For simplicity, we consider only the case of a single \app{} of two
parties; larger \app{s}  can be treated similarly as well as several \app{s}  
(their effects are, asymptotically at least, additive).

Let the parties be $i$ and $j$, and consider the expectation of the gain
$\sij-(s_i+s_j)$, where $\sij$ is the 
number of 
seats the parties get together as an \app.

\begin{theorem}\label{Tmerge}
  Suppose that $\ppm$ are \linQ. 
If two parties $i$ and $j$ form an \app{} in the election,
the mean of their seat gain $\sij-s_i-s_j$ satisfies, as $\Ntoox$,
\begin{align}\label{tmd}
  \E(\sij-s_i-s_j)
&\to (\gb-\tfrac12)\bigpar{1-p_i-p_j}
\intertext{for the \gbdiv{}, and}
\label{tmq}
  \E(\sij-s_i-s_j)
&\to \gam\Bigpar{\frac{2}{m}-\frac{1}{m-1}}
= \gam\frac{m-2}{m(m-1)}
\end{align}
for the \gamq{} method.
\end{theorem}

\begin{proof}
  If $\gdij$ is the seat excess for the  \app{}, then
  \begin{equation*}
	\sij-s_i-s_j
=N(p_i+p_j)+\gdij-(Np_i+\gD_i)-(Np_j+\gD_j)
=\gdij-\gD_i-\gD_j,
  \end{equation*}
so 
$  \E(\sij-s_i-s_j)=\E\gdij-\E\gD_i-\E\gD_j$, which
by Theorems \refand{TDmean}{TQmean}
converges to,
recalling that $\gdij$ is the seat excess in a contest with   $m-1$ parties,
\begin{equation*}
(\gb-\tfrac12)\Bigpar{(m-1)(p_i+p_j)-1-(mp_i-1)-(mp_j-1)}
\end{equation*}
and
\begin{equation*}
\gam\Bigpar{\Bigpar{p_i+p_j-\frac{1}{m-1}}
-\Bigpar{p_i-\frac{1}{m}}
-\Bigpar{p_j-\frac{1}{m}}}
\end{equation*}
for the \gbdiv{} and \gamq{} method, respectively.
\end{proof}
One sees in the same way that these expected gains (which are negative if
$\gb<1/2$ or $\gam<0$) are balanced by a loss for each other party $k$
of $(\gb-\frac12)p_k$ for the \gbdiv{} and $\gam/(m(m-1))$ for the \gamq{}
method. 

We have seen that \app{s} are favoured when $\gb>1/2$ or $\gam>0$, but
the effect is rather small, especially for quota methods.
(Nevertheless, as said above, 
the effect is large enough to be politically important 
in reality for D'Hondt's method.)
Note that for the \gbdiv, the gain in \eqref{tmd} does not depend on
the number (or sizes) ot the other parties, while for the \gamqm, the gain
in \eqref{tmq} is largest for $m=3$ or 4 and decreases as $O(1/m)$ for large
$m$.

\begin{example}
  For Jefferson/D'Hondt's method, two small parties that form an \app{}
gain at most 0.5 seats together, and two large parties gain less.
\end{example}

\begin{example}
  For Droop's method,
any two parties that form an \app{}
gain at most $\frac16\approx0.167$ seats together when there are 3 or 4
parties,  $0.15$ seats when there are 5 parties, and less if there are more.
\end{example}

\begin{problem}
  What are, for the different methods,  the asymptotic distributions
of the seat gain $\sij-s_i-s_j$ for  an \app{}? 
\end{problem}

We have so far studied the gain of the two parties combined, but the parties
are probably more interested in how the gain is split between them.
Typically, the seats given to an \app{} are distributed,
in a \emph{sub-apportionment}, between the
participating parties by the same election method as for the main
distribution (the \emph{super-apportionment}).
It seems reasonable that the gain then, on the average, is split between the
parties proportionally to their sizes, \cf{} Theorems \refand{TDND}{TDNQ}.
However, we shall see that this holds only for the divisor methods.

Consider first the \gbdiv{}. 
\refT{TDmean} shows that, asymptotically as \Ntoox, the number of seats the
\app{} gets is on the average
\begin{equation}
  N(p_i+p_j)+(\gb-\tfrac12)\bigpar{(m-1)(p_i+p_j)-1}.
\end{equation}
If we use \refT{TDmean} again for the distribution of these seats between
$i$ and $j$, we obtain that $i$ gets on the average
{\multlinegap=0pt%
\begin{multline}
\Bigpar{ N(p_i+p_j)+(\gb-\tfrac12)\bigpar{(m-1)(p_i+p_j)-1}}\frac{p_i}{p_i+p_j}
+(\gb-\tfrac12)\Bigpar{2 \frac{p_i}{p_i+p_j}-1}
\\
=
Np_i+  
(\gb-\tfrac12)\Bigpar{(m-1)p_i+\frac{p_i}{p_i+p_j}-1}
\end{multline}}%
seats, which, comparing to \refT{TDmean} for an election without the \app,
means an average gain for party $i$ of 
\begin{equation}\label{tm}
(\gb-\tfrac12)\Bigpar{\frac{p_i}{p_i+p_j}-p_i}
=
\frac{p_i}{p_i+p_j}(\gb-\tfrac12)\xpar{1-p_i-p_j},
\end{equation}
which indeed is the proportional share $p_i/(p_i+p_j)$ of the joint gain in
\eqref{tmd}. 

There is, however, a gap in this argument. In our model, the seat assignment
is a deterministic function of $N$ both in the super-apportionment and in
the sub-apportionment, and the calculation above using expectations
for random $N$ assumes that there is no hidden correlation, where
the numbers of seats given to the \app{} are seat
numbers that tend to favour one of the two parties in the
sub-apportionment. It seems very unlikely that there will be such a
correlation under our assumption that $\ppm$ are \linQ, but we have so far
no rigorous proof, so \eqref{tm} should only be regarded as a conjecture for
the asymptotic average gain for party $i$.

\begin{problem}
  Verify rigorously that the expected gain for party $i$ converges to the
  value in \eqref{tm} as \Ntoox.
\end{problem}

Nevertheless, if we tentatively use \eqref{tm},
and regard it as a valid approximation for real elections with finite $N$
(and arbitrary $\ppm$), we can draw some practical conclusions for
D'Hondt's method ($\gb=1$).
 If the parties $i$ and $j$ both are small, then \eqref{tm} shows that the
 gain for party $i$ is almost $p_i/(p_i+p_j)$. Thus most of the gain of the
 \app{} goes to the larger party, and a party has very little to gain by an
 \app{} with a party that is substantially larger.
In fact, if party $i$ can choose between several parties to form an \app, we
see, perhaps surprisingly, that it gets the largest gain by choosing the
smallest partner. (This is based on asymptotics as $\Ntoo$. For a real
election with finite $N$ there is a trade-off since a really small partner
will not help significantly; the advantage can hardly be larger than $Np_j$,
so presumably only parties with $q_j=Np_j$ being at least $1/2$ or so will
be useful partners. It would be interesting to study such effects for finite
$N$, but that is beyond the scope of the present paper.)

\smallskip
Consider now instead the \gamqm. 
One important complication is that quota methods are not uniform, in the
sense of \refApp{Auniform}. Thus, even if the \app{} gets the same number of
seats as the two parties would have got together if they had appeared
separately in the election, the sub-apportionment may give a different
distribution of these seats between the parties.

We calculate as above, with the same reservation that a rigorous
verification still is lacking.
\refT{TQmean} shows that, asymptotically as \Ntoox, the number of seats the
\app{} gets is on the average
\begin{equation}
  N(p_i+p_j)+\gam\Bigpar{(p_i+p_j)-\frac1{m-1}}.
\end{equation}
If we use \refT{TQmean} again for the distribution of these seats between
$i$ and $j$, we obtain that $i$ gets on the average
{\multlinegap=0pt%
\begin{multline}
\Bigpar{  N(p_i+p_j)+\gam\Bigpar{(p_i+p_j)-\frac1{m-1}}}\frac{p_i}{p_i+p_j}
+\gam\Bigpar{ \frac{p_i}{p_i+p_j}-\frac12}
\\
=
Np_i+  
\gam\Bigpar{p_i+\frac{m-2}{m-1}\frac{p_i}{p_i+p_j}-\frac12}
\end{multline}}%
seats, which, compared to \refT{TQmean} for an election without the \app,
means an average gain for party $i$ of 
\begin{equation}\label{tmq2}
\gam\Bigpar{\frac{m-2}{m-1}\frac{p_i}{p_i+p_j}-\frac{m-2}{2m}}.
\end{equation}
This is quite different from the proportional share $p_i/(p_i+p_j)$ of the
joint gain in \eqref{tmd}. 
In fact, we see that for $\gam>0$ (\eg{} Droop's method),
a party gains only if $p_i/(p_i+p_j)>(m-1)/(2m)$.
In other words, a party will \emph{lose} if it forms an \app{} with a party
that is only a little larger. Typically, the redistribution within the
\app{}
by the sub-apportionment is more important than the collective gain.

\section{The \StL{} divergence}\label{Sdiv}

\citet{StL1,StL2} based his proposal of his method on the fact that it
is the least squares method
minimising the sum
\begin{equation}
S\=\sumim v_i\Bigpar{\frac{s_i}{v_i}-\frac NV}^2
=
\sumim \frac{(s_i-Np_i)^2}{v_i}
=
\sumim \frac{\gD_i^2}{v_i}
;  
\end{equation}
multiplying this by $V$ we obtain the
equivalent quantity
\begin{equation}\label{sld}
  \sld\=VS=\sumim\frac{\gD_i^2}{p_i},
\end{equation}
which is called the \emph{\StL{} divergence} by \citet{HPS2004,HPS2005},
where the 
asymptotic distribution is studied under the assumption of random party
sizes
(see \refS{Srandom}). We obtain corresponding result in our setting as
corollaries of the results in \refS{Sresults}.
(Unfortunately, we do not get as explicit result for quota methods as for
divisor methods because \refT{TQjoint} is less explicit.)
The most interesting case is the $\frac12$-linear divisor method (\StL's
method)  \cite{HPS2004}, 
since then $S$ and $\sld$ are the minima over all allocations,
but it is also interesting to see how much larger $\sld$ is for other
methods \cite{HPS2005}.

\begin{theorem}\label{Tsld}
  Suppose that $\ppm$ are \linQ, and let \Ntoox.
  \begin{romenumerate}
  \item 
For the \gbdiv,
\begin{equation}\label{tsldd}
  \sld \dto
\sbm\=
\sumim\frac{(V_i+\gb-1)^2}{p_i}-\Bigpar{\sumim(V_i+\gb-1)}^2,
\end{equation}
where $V_1,\dots,V_m$ are as in \refT{TD}.
\item 
For the \gamqm,
\begin{equation}\label{tsldq}
  \sld \dto
\sgm\=
\sumim\frac{\bY_i^2}{p_i},
\end{equation}
where $\bY_1,\dots,\bY_m$ are as in \refT{TQjoint}.
  \end{romenumerate}
\end{theorem}

\begin{proof}
  For the \gbdiv, we use \refT{TD} and obtain 
by the continuous mapping theorem
\cite[Theorem 5.1]{Billingsley} 
and \eqref{sld} immediately
$\sld\dto\sumim X_i^2/p_i$.
We write \eqref{tdxi} as $X_i=p_i A-A_i$, where 
$A_i\=V_i+\gb-1$ and 
$A\=\sumim (V_i+\gb-1)=\sumim A_i$.
Thus
\begin{equation*}
  \begin{split}
	\sumim\frac{X_i^2}{p_i}
&=	\sumim\frac{p_i^2A^2-2p_iAA_i+A_i^2}{p_i}
=	\sumim p_i A^2-2\sumim A_i A+\sumim\frac{A_i^2}{p_i}
\\&
=\sumim\frac{A_i^2}{p_i}-A^2,
  \end{split}
\end{equation*}
which is \eqref{tsldd}. 

The result for the \gamqm{} follows directly from \refT{TQjoint}.
\end{proof}

The expression in \eqref{tsldd} for $\sbm$ may be compared to the somewhat
different expression in \cite{HPS2004,HPS2005} for the limit for fixed $N$
and random $\ppm$; by the argument in \refS{Srandom}, the formulas
have to be equivalent (for random $\ppm$).

For the expectation we similarly obtain:

\begin{theorem}\label{Tsldmean}
  Suppose that $\ppm$ are \linQ, and let \Ntoox.
  \begin{romenumerate}
  \item 
For the \gbdiv,
\begin{equation}\label{tslddmean}
\E  \sld \to
\frac1{12}\lrpar{\sumim\frac{1}{p_i}+m-2}
 +(\gb-\tfrac12)^2\lrpar{\sumim\frac1{p_i}-m^2}.
\end{equation}
\item 
For the \gamqm,
\begin{equation}\label{tsldqmean}
\E  \sld 
\to
\frac{(m+2)(m-1)}{12 m^2}\sumim\frac{1}{p_i}
 +\frac{\gam^2}{m^2}\lrpar{\sumim\frac1{p_i}-m^2}.
\end{equation}
  \end{romenumerate}
\end{theorem}

\begin{proof}
  This can be shown from \refT{Tsld}, but it is easier to use
\begin{equation*}
\E  \sld=\sumim\frac{\E\gD_i^2}{p_i}
=\sumim\frac{\Var\gD_i+(\E\gD_i)^2}{p_i}
\end{equation*}
and \refT{TDmean}, \refC{CDvar}, \refT{TQmean}, \refC{CQvar}.
\end{proof}

Note that $\sumim 1/p_i\ge m^2$, with equality if and only if all $p_i$ are
equal. 
We know that $\sld$ always is minimised by the $\frac12$-linear divisor
method (\StL's method).
We see from \refT{Tsldmean} that also asymptotically, $\E\sld$ is strictly
larger for any \gbdiv{} or \gamqm, 
except when all parties have exactly the same size (a trivial case where all
methods yield the same result; note that we
have excluded this case by our assumptions).

Furthermore, Theorems \ref{Tsld}--\ref{Tsldmean} show that
$\sld$ is typically of the order $m^2$ if all $p_i$ are of
roughly the same order, but if some $p_i$ is very small, $\sld$ can be
substantially larger.
\citet{HPS2004,HPS2005} consider $\sld$ and the limit $\sbm$
for random, uniformly distributed,
$\ppm$, and show that, as \mtoo,
\begin{equation}\label{1stable}
\xfrac{\sbm}{m^2}  -\bigpar{\tfrac1{12}+(\gb-\tfrac12)^2}\log m \dto \sldlim,
\end{equation}
where $\sldlim$ is a certain 1-stable random variable (depending on $\gb$).
We see that the $\log m$ term and the large tail of $\sldlim$ 
(we have $\E\sldlim=\infty$)
depend on the
possibility of some very small $p_i$; the asymptotic behaviour
for deterministic $\ppm$ with, say, $(\min p_i )\qw = O(m)$ is quite
different.
To see this better, we first replace $V_i$ (which are dependent) by $U_i$
(which are \iid) in \eqref{tsldd}; the following lemma shows that for large
$m$, the difference is small.
(Cf.\ 
similar approximations in 
\cite{HPS2004,HPS2005}.)

\begin{lemma}
  Let $\gbu_i\=U_i+\gb-1$, 
so $\gbu_i$ are \iid{}
  $\U(\gb-1,\gb)$, and let 
  \begin{equation}\label{thk}
	\sbm'\=
\sumim\frac{{\gbu_i}^2}{p_i}- \biggpar{\sumim {\gbu_i}}^2.
  \end{equation}
Then
$\E\bigabs{\sbm-\sbm'}=O(m)$, uniformly in $\ppm$.
\end{lemma}
\begin{proof}
  We have $|U_i+\gb-1|\le|\gb|+1=O(1)$ and similarly
$|V_i+\gb-1|=O(1)$. Further, by \eqref{tdvi}, $\P(V_i\neq  U_i)=\P(J=i)=p_i$.
Hence,
\begin{equation}\label{thk1}
  \E\frac{\lrabs{(V_i+\gb-1)^2-(U_i+\gb-1)^2}}{p_i}
=\frac{O(p_i)}{p_i}=O(1).
\end{equation}
Similarly,
$\sumim (V_i+\gb-1)-\sumim (U_i+\gb-1)=-U_J=O(1)$ and thus
\begin{equation}\label{thk2}
\biggpar{\sumim (V_i+\gb-1)}^2-\biggpar{\sumim (U_i+\gb-1)}^2
=O(m).  
\end{equation}
The result follows by summing \eqref{thk1} and \eqref{thk2}.\end{proof}

We can now state a normal limit theorem corresponding to \eqref{1stable},
but for deterministic $\ppm$.

\begin{theorem}
  Suppose that for each $m$, we are given $\ppm$ such that
  $m^{3/2}\min_ip_i\to\infty$ as \mtoo. 
Let 
$A_m\=\sumim p_i\qw$,
$B_m\=\sumim p_i\qww$ and
\begin{align*}
C_m&\=\sumim (p_i\qw-m)^2=B_m-2mA_m+m^3.
\end{align*}
Then, as \mtoo,
\begin{equation*}
  \frac{\sbm-\bigpar{A_m/12+(A_m-m^2)(\gb-1/2)^2}}
{\sqrt{B_m/180+(\gb-1/2)^2C_m/3}}
\dto
N\bigpar{0,1}.
\end{equation*}
\end{theorem}

\begin{proof}
  Note first that, \eg{} by Jensen's inequality, 
$A_m\ge m^2$, $B_m\ge m^3$ and $C_m\ge0$.

Write, for convenience, $b\=\gb-1/2$.
We have $\gbu_i=\tU_i+b$ and $\E\gbu=b$.
Thus, \eqref{thk} can be expanded as
\begin{equation*}
  \begin{split}
		\sbm'&=
\sumim\frac{\tU_i^2+2b\tU_i+b^2}{p_i}- \biggpar{mb+\sumim {\tU_i}}^2
\\&=
\sumim\Bigpar{\frac{\tU_i^2+2b\tU_i}{p_i}-2mb\tU_i}
+b^2A_m-b^2m^2- \biggpar{\sumim {\tU_i}}^2.
  \end{split}
  \end{equation*}
The last term has expectation $\Var\sumim \tU_i=O(m)$, and can thus be ignored. 
Define
\begin{equation*}
  \XX_i\=\frac{\tU_i^2+2b\tU_i}{p_i}-2mb\tU_i.
\end{equation*}
Then 
\begin{equation*}
\E  \XX_i=\frac{\E\tU_i^2}{p_i}=\frac{1}{12p_i}
\end{equation*}
and, by a straightforward calculation,
\begin{equation*}
\Var \XX_i=
\frac1{180p_i^2}+\frac{b^2}3\Bigpar{\frac1{p_i}-m}^2.
\end{equation*}
The random variables $\XX_i$ are independent and each is bounded with
$|\XX_i|=O(1/p_i+m)=o(m\qqc)=o(B_m\qq)$.
The result now follows by the central limit theorem, 
using either Lindeberg's or Lyapounov's condition,
see \eg{} \cite[Section 7.2]{Gut}.
\end{proof}

\begin{remark}
\citet{StL2} also studied the least squares functional calculated per seat
and not per voter:
\begin{equation}
  \sumim s_i\Bigpar{\frac{v_i}{s_i}-\frac VN}^2
=\frac{V^2}{N^2}\sumim\frac{\gD_i^2}{s_i},
\end{equation}
and found that this is minimised by the divisor method with
$d(n)=\sqrt{n(n-1)}$, later suggested by Huntington, see \refApp{Adivisor}
and \refR{RHuntington}.
The asymptotics of this quantity is, after division by $(V/N)^3$, the same
as for $S$, because $s_i/N\to p_i$.
\end{remark}

\section{Some other goodness-of-fit functionals}\label{Sothergood}
Many election methods, including the ones treated here, 
can be characterised as minimising various
functionals, see \eg{} \cite{BY}, \cite{Kopfermann} and \cite{NN}.
Theorems \refand{TD}{TQjoint} yield by the continuous mapping theorem
asymptotic distributions for many such functionals for the methods treated
here, under our standard assumptions that $\ppm$ are \linQ{} and $\Ntoox$,
which we assume below.
We give some examples,
leaving others to the reader.
Note, however, that our results for quota methods are
less complete, \cf{} \refP{PQ}.

\subsection{Hamilton divergences}

Hamilton/Hare's method, the $0$-quota me\-thod, is the unique method
minimising 
$\max_i|\gD_i|$. It also the method minimising $\max_i\gD_i$.
Moreover, it is,
for any convex function $\gf:\bbR\to[0,\infty]$ with $\gf(0)=0$ and 
$\gf(x)>0$ for $x\neq0$, 
the unique method minimising
$\sum_i\gf(\gD_i)$, 
as was shown by \citet{Polya1918b}. 

For the \gbdiv, \refT{TD} implies
\begin{equation}
  \maxim|\gD_i|\dto \maxim |X_i|
\end{equation}
with $X_i$ given by \eqref{tdxi}. However, we do not know any really simple
description of this limit variable.

For the \gamqm, including Hamilton/Hare's method and thus the minimum of
$\max_i|\gD_i|$ over all allocations, \refT{TQjoint} yields 
\begin{equation}
  \maxim|\gD_i|\dto \maxim |Y_i|,
\end{equation}
but we have no explicit description of this limit  variable.

Similar results hold for $\max_i\gD_i$.

We obtain more complete results if we instead consider the least squares
functional $\sumim\gD_i^2$, which also yields Hamilton/Hare's method by
Polya's theorem \cite{Polya1918b}. (This case was shown earlier by
\citet{StL2}, who attributed the result to
Zivy.)

\begin{theorem}
  With assumptions and notations as in \refS{Sresults}, for the \gbdiv,
  \begin{align}
\sumim\gD_i^2 &\dto\sbm\hh\=\sumim X_i^2,\\	
\E\sumim\gD_i^2 &\to\E\sbm\hh
=\sumim\frac{1+(m-2)p_i^2}{12}+(\gb-\frac12)^2\sumim(mp_i-1)^2
\\
&
=\frac{(m+2)(m-1)}{12m}+\Bigpar{\frac{m-2}{12}+m^2(\gb-\frac12)^2}
\sumim\Bigpar{p_i-\frac1m}^2,
  \end{align}
and for the \gamqm,
  \begin{align}
\sumim\gD_i^2 &\dto\sgm\hh\=\sumim \bY_i^2,\\	
\E\sumim\gD_i^2 &\to\E\sgm\hh
=\frac{(m+2)(m-1)}{12m}+\gam^2\sumim\Bigpar{p_i-\frac1m}^2.
  \end{align}
\end{theorem}

\begin{proof}
By Theorems \refand{TD}{TQjoint}, calculating $\E X_i^2$ and $\E \bY_i^2$ by
 \refT{TDmean}, \refC{CDvar}, \refT{TQmean}, \refC{CQvar} as in
the proof of \refT{Tsldmean}.
\end{proof}

\subsection{Jefferson and Adams divergences}

\citet{HS} consider the functionals
\begin{align}
  S\jj&\=\maxim\Bigpar{\frac{s_i}{p_i}-N}=\maxim\frac{\gD_i}{p_i}
\intertext{and}
  S\qa&\=\minim\Bigpar{\frac{s_i}{p_i}-N}=\minim\frac{\gD_i}{p_i}.
\end{align}
(We could as well use the maximum and minimum of $s_i/v_i-N/V$.)
Since $\sumim \gD_i=0$, we have $S\jj\ge0$ but $S\qa\le0$, and we will
therefore  use $|S\qa|=-S\qa$.

It is easy to see that $\sjj$ and $\saa$ are minimised, over all allocations
of $N$ seats, by the Jefferson and Adams methods, respectively \cite{StL2},
\cite{BY}. 

\refT{TD} yields the asymptotic distribution of $\sjj$ and $\saa$ for any
\gbdiv; we consider for simplicity only the optimal methods.

\begin{theorem}
  \label{TJA}
With assumptions and notations as in \refS{Sresults}, 
for the optimal allocations,
\begin{align}
  \sjj&\dto\sjjm\=\sumjmi U_j,
\label{tjaj}
\intertext{and}
  \saa&\dto\saam\=m-\sumjm V_j -\minim\frac{1-V_i}{p_i}.
\label{tjaa}
\end{align}
We have
\begin{equation}
  \label{tje}
\E\sjjm=\E\saam=\frac{m-1}2.
\end{equation}
\end{theorem}

\begin{proof}
  $\sjj$ is minimised by Jefferson's method, so \refT{TD} with $\gb=1$
  yields, noting that $\min V_i/p_i=0$ by \eqref{tdvi},
  \begin{equation*}
	\begin{split}
\sjj&\=\maxim\frac{\gD_i}{p_i}
\dto\maxim\frac{X_i}{p_i}
=\sumjm V_j-\minim \frac{V_i}{p_i}
\\&\phantom:
=\sumjm V_j
=\sum_{j\neq J} U_j
\eqd\sumjmi U_j.
	\end{split}
  \end{equation*}

Similarly, $\saa$ is minimised by Adams's method, so \refT{TD} with $\gb=0$
yields
  \begin{equation*}
	\begin{split}
\saa&\=-\minim\frac{\gD_i}{p_i}
\dto-\minim\frac{X_i}{p_i}
=-\sumjm V_j+m-\minim \frac{1-V_i}{p_i}.
	\end{split}
  \end{equation*}

We have $\E\sjjm=(m-1)/2$ directly from \eqref{tjaj}.

For $\saam$, we first note that 
if $B\=\minim 1/p_i$, then for $0<x<B$,
\begin{equation*}
  \begin{split}
\P\Bigpar{\minim \frac{1-V_i}{p_i}>x}
&=\P\bigpar{V_i<1-xp_i,\; i=1,\dots,m}	
\\&
=\sumjm p_j\P\bigpar{V_i<1-xp_i,\; i=1,\dots,m\mid J=j}	
\\&
=\sumjm p_j\P\bigpar{U_i<1-xp_i,\; i\neq j}
\\&
=\sumjm p_j\prod_{i\neq j}(1-xp_i),
  \end{split}
\end{equation*}
while $\P\Bigpar{\minim \xfrac{(1-V_i)}{p_i}>x}=0$ for $x\ge B$.
Hence,
\begin{multline}\label{eva}
\E\minim \frac{1-V_i}{p_i}
=\intoo \P\Bigpar{\minim \frac{1-V_i}{p_i}>x}\dd x
\\
=\int_0^B \sumjm p_j\prod_{i\neq j}(1-xp_i) \dd x
=\lrsqpar{-\prod_{i=1}^m(1-xp_i)}_{x=0}^B
=1.
\end{multline}

Furthermore, $\E\sumjm V_j=\sumjmi U_j=(m-1)/2$, and \eqref{tje} follows
from \eqref{tjaa}.
\end{proof}

\begin{remark}\label{RSA}
  Note that $\sjjm$ in \eqref{tjaj} does not depend on $\ppm$, but for
	$\saam$, the situation is more complicated; the representation in
	\eqref{tjaa} depends on $\ppm$, but it is still possible that $\saam$
	has the same distribution for all $\ppm$ and that we just have found an
 unnecessarily complicated description of it.
Note that $\E\saam$ does not depend on $\ppm$. In the case $m=2$, it is easy
to see that $\saam\eqd U\eqd\sjjm$ for every $(p_1,p_2)$, but we do not know
whether something similar is true for $m\ge3$.
\end{remark}

\begin{problem}
  Is $\saam\eqd\sjjm$ for every $\ppm$?
\end{problem}

The asymptotic results by \cite{HS}, taking the limit of the limits as
\mtoo, follow easily in our setting too.

\begin{theorem}
  For any $\ppm$ depending on $m$, as \mtoo,
  \begin{align*}
	\sqrt m(\sjjm-m/2)&\dto N(0,\tfrac1{12}),
\\
	\sqrt m(\saam-m/2)&\dto N(0,\tfrac1{12}).
  \end{align*}
\end{theorem}

\begin{proof}
  The result for $\sjjm$ is immediate from \eqref{tjaj} and the central
  limit theorem.

For $\saam$ we use \eqref{tjaa}. We have
${\sumjm(U_j-V_j)}=U_J$ and
consequently, by \eqref{tjaa} and \eqref{eva},
\begin{equation*}
  \begin{split}
\E\biggabs{\saam-\sumjm(1-U_j)}
&\le \E \Biggabs{\sumjm(U_j-V_j)-\minim \frac{1-V_i}{p_i}}
\\&
\le \E U_J + \E\minim \frac{1-V_i}{p_i}
<2,	
  \end{split}
\end{equation*}
and the result follows by the central limit theorem applied to $\sumjm(1-U_j)$.
\end{proof}

\section{Rational $p_i$}\label{Sdep}

In our main results we have assumed that $\ppm$ are \linQ, since this is
necessary for \refL{LW}.
We consider briefly what happens when this is not satisfied, \ie, when
$\ppm$ are linearly dependent over $\bbQ$. In particular, we are for obvious
practical reasons interested in the case when $\ppm$ are rational.
(See also the corresponding discussion in \cite{SJ258}.)

Note first that
it may now happen that there are ties. We assume that any ties are resolved
randomly; thus the quantities $s_i(N)$ and $\gD_i(N)$ in general may be
random variables, also for a fixed $N$. This is no real problem, however,
and when discussing the means we may replace them by their average over all
solutions in case of a tie.

\refL{LW} does not apply when
$p_1,\dots,p_m$ are linearly dependent over $\bbQ$, but the proof of it sketched
in \refS{SpfD} shows that the sequence 
$(\frax{ny_1+a_1},  \dots,\frax{ny_k+a_k})_{n\ge 1}$ 
is uniformly distributed on a coset of a subgroup
of $[0,1)^{k}$; more precisely, if we for simplicity take all $a_i=0$,
the empirical distributions converge to 
the uniform probability measure $\mu_\yy$ on a
subgroup, with Fourier coefficients given by
\begin{equation}\label{mupp}
\widehat\mu_\yy(\ell_1,\dots,\ell_k)=
\begin{cases}
  1 & \text{ if } \sum_{j=1}^k\ell_jy_j\in\bbZ,\\
0 & \text{ otherwise}.
\end{cases}
\end{equation}
The proofs of Theorems \refand{TD}{TQjoint} now show, with minor
modifications, that the seat excesses
$\gD_i$ still converge jointly in distribution as $\Ntoox$, but the limit
distribution is different, and depends on $\ppm$ through some measures of
the type \eqref{mupp}.
In particular, the bias $\lim\E\gD_i$ is well-defined, but we do not have an
explicit formula for it. 

However, if we consider a
sequence of party size distributions
$(p_{1k},\dots,p_{mk})$, $k=1,2,\dots$, (with a fixed $m$), such that
 $p_{ik}\to p_i$ for each $i$ as \ktoo{}
and further,
for every integer vector $(a_1,\dots,a_m)\neq0$, we have 
$\sum_{i=1}^m a_ip_{ik}\neq0$ for all large $k$,
then it follows from \eqref{mupp} and an inspection of
the proofs of Theorems \refand{TD}{TQjoint} 
that the measures $\mu_\yy$ that now appear in the proofs will converge
to the uniform measure $\mu$, and thus the limit variables $X_i$ and $\bY_i$
will converge to the corresponding variables given in \refS{Sresults} for
the case when $\ppm$ are \linQ.

This shows that the results above are good approximations also for linearly
dependent $\ppm$, as long as there is no
linear relation \eqref{linq} with small integers $a_i$.
In particular, for rational $\ppm$ we typically have a good approximation
unless some $p_i$ has a small denominator.

In the remainder of this section, we suppose that $\ppm$ are rational,
with a common denominator $L$.
Theorems \refand{TDND}{TDNQ} then show that the sequence of seat excesses 
$\gD_i\=s_i(N)-Np_i$, $N\ge1$, has period $L$. 
(If $\gb<0$, $\gb>1$ or $\gam<0$ we may have to except some small $N$.)
This gives another (simpler)
proof that 
the limit distribution of $\gD_i$ as $\Ntoox$ exists in this case, and shows
that it is the same as the 
distribution if we take $N$ uniformly in 
a period \set{K+1,\dots,K+L} (for any $K$ that is large enough).
In particular, the
asymptotic distribution is discrete, and therefore not the same as in
\refS{Sresults}. The asymptotic 
bias $\lim \E\gD_i$ is easily calculated by taking the
average of $\gD_i$ over a period.

\begin{example}\label{E13}
  Let $m=2$ and take $(p_1,p_2)=(\frac23,\frac13)$. 
Then the \gbdiv{} equals the \gamqm{} with $\gam=2\gb-1$, see
  \refApp{SS2}. 
Assume, for simplicity, $0<\gb< 1$, \ie, $-1< \gam < 1$.
A simple  calculation shows that for $N=1,2,3$, we have
$\gD_1=\frac13,-\frac13,0$. Hence the average is 0, so the methods are
unbiased in this case, for any $\gb\in(0,1)$ or $\gam\in(-1,1)$. (In fact,
for this example the methods all coincide with Hamilton/Hare =
Webster/\StL{}.) In particular, Theorems \refand{TDmean}{TQmean} do not hold
for rational \ppm{} and $\gb\in(0,\frac12)\cup(\frac12,1)$ and
$\gam\in(-1,0)\cup(0,1)$. 

If we instead take $\gb=1$, \ie{} $\gam=1$, we have a tie for $N=2$, with
$\gD_1=-\frac13$ or $\frac23$. The average of $\gD_1$ over a period now is
$\frac16$.
\end{example}

For the methods with $\gb=1/2$ and $\gam=0$, the result 0 in \refE{E13}
agrees with Theorems \refand{TDmean}{TQmean}. In fact, this is true for any
rational $\ppm$.

\begin{theorem}\label{TSW}
  The $\frac12$-linear divisor method (Webster/\StL) and the $0$-quota  method 
(largest remainder/Hamilton/Hare)
are asymptotically unbiased for any rational $\ppm$.
\end{theorem}
\begin{proof}
  Let, as above, $L$ be a common denominator of $\ppm$. It is for both
  methods easily seen, arguing similarly to the proofs of Theorems
  \refand{TDND}{TDNQ}, 
  that $s_i(L-N)=Lp_i-s_i(N)$ for $N=0,\dots,L$. Hence,
$\gD_i(L-N)=-\gD_i(N)$, and the average over a period $N=1,\dots,L$ vanishes.
\end{proof}

Also for $\gb=1$ and $\gam=1$,  the result in \refE{E13}
agrees with Theorems \refand{TDmean}{TQmean}. 
For divisor methods, we have the following general result.
(By \refE{E13} and \eqref{rgbper}, the result does not hold for any other
$\gb$.) 

\begin{theorem}
  \label{TErika}
Consider the \gbdiv{} and suppose that $\gb=k/2$ for some integer $k$. Then
\eqref{tdmean} holds for any rational \ppm.
\end{theorem}
\begin{proof}
  By \eqref{rgbper}, the result holds for $\gb$ if and only if it holds for
  $\gb+1$. Hence it suffices to consider $\gb=1/2$ and $\gb=1$.
The case $\gb=1/2$ is part of \refT{TSW}.

For $\gb=1$, we argue similarly. (Cf.\ the proof of \refT{TDND}.)
Let $0\le N<L$ and let $t\ge0$ be such that $s_i(N)=\rgx1{p_it}$, see
\eqref{gbdiv} and \eqref{divam}. Recall that $\rgx1{p_it}=\floor{p_it}$
except that $\floor{p_it}-1$ also is possible when $p_it$ is an integer.

Since $p_iL$ is an integer, if we first for simplicity assume that $p_it$ is
not, 
\begin{equation}\label{th}
\rgx1{p_i(2L-t)}=2p_iL-\rgx1{p_it}-1=2p_iL-s_i(N)-1,  
\end{equation}
where summing over $i$ yields
\begin{equation}
\sumim\rgx1{p_i(2L-t)}=2L-\sumim s_i(N)-m=2L-N-m.
\end{equation}
Consequently, \eqref{divam} again 
shows that these numbers yield the seat distribution for $2L-N-m$ seats.
(Note that $L=\sumim p_iL\ge m$, and thus $2L-N-m>0$.)
Hence, using \eqref{th} again,
\begin{equation}\label{kk}
s_i(2L-m-N)=
\rgx1{p_i(2L-t)}
=2p_iL-s_i(N)-1.
\end{equation}
The argument works also, with a little care, 
if some $p_it$ is an integer, and shows that
\eqref{kk} holds generally, and that ties in the distribution of $s_i(N)$  
correspond to ties in the distribution of
$s_i(2L-m-N)$.  
Hence, we have, regarding $\gD_i$ as a random variable when there is a tie,
\begin{equation}
\gD_i(2L-m-N)=
s_i(2L-m-N)-(2L-m-N)p_i
=mp_i-\gD_i(N)-1.
\end{equation}
Taking the expectations in case of a tie, and 
then the average over the period $N=0,\dots,L-1$, when $2L-m-N$ runs
through the period $L-m+1,\dots,2L-m$, we see that if the average is $x$,
then $x=mp_i-x-1$; hence $x=\frac12(mp_i-1)$, and thus \eqref{tdmean}
holds when $\gb=1$.
\end{proof}

For the 1-quota method (Droop), we have seen that \eqref{tqmean} holds 
in \refE{E13}; in fact, when $m=2$, it holds for any 
$(p_1,p_2)$ by \refT{TErika}, because the 1-quota method equals the 1-linear
divisor method when $m=2$, see \refApp{SS2}.
However, it does not hold for $m=3$.

\begin{example}
  \label{E122}
Let $m=3$ and $(p_1,p_2,p_3)=(\frac25,\frac25,\frac15)$. 
Consider the 1-quota method (Droop).
For $N=1,2,3,4,5$, the smallest party gets,
taking the average when there is a tie,
 $0,0,1,\frac23,1$ seats. (Note the non-monoticity, an instance of the
Alabama paradox which frequently occurs for all quota methods, \cf{}
  \cite{SJ258}.) 
The seat excesses are $-\frac15,-\frac25,+\frac25,-\frac2{15},0$, with
average $-\frac1{15}$, while \eqref{tqmean} would give 
$p_1-\frac13=-\frac2{15}$.
Thus \eqref{tqmean} does not hold.

Note that this example shows that the heuristic argument in \refR{RQheur} is
not always valid when \ppm{} are rational. In this example, the retracted
seat in \refR{RQheur} is \emph{not} taken uniformly; a calculation shows
that it is taken from the smallest party only $\frac4{15}$ of the times.
\end{example}

We can ask the following for any \gbdiv{} or \gamqm,
where the cases $\gam=1$ in \refP{P1} and 
$\gb=1/2$, $\gb=1$ and
$\gam=0$ in \refP{P2}
are especially interesting.
(By \refT{Tgb+1}, for the linear divisor
methods, it suffices to consider $0<\gb\le1$.)

\begin{problem}\label{P1}
  Find a general formula for the asymptotic bias $\lim\E\gD_i$ for rational
  \ppm. 
(Theorems \refand{TSW}{TErika} yield the answers for $\gb\in\bbZ/2$ and
  $\gam=0$.) 
\end{problem}

\begin{problem}\label{P2}
  Find a general formula for the asymptotic distribution of $\gD_i$ for rational
  \ppm. 
\end{problem}

\appendix

\section{Divisor methods}\label{Adivisor}

A \emph{divisor method} is based on a given sequence of numbers
$d(n)$, $n\ge1$, with $0\le d(1)<d(2)<d(3)<\dots$.
Different choices of $d(n)$ give different methods.
A number of divisor methods that have been used or discussed are shown in
\refTab{TabDiv}. The most important ones are the widely used methods by
Jefferson/D'Hondt and Webster/\StL. 
(Notation varies. In \eg{} \cite{BY}, the sequence is denoted
$d(0),d(1),\dots$; thus their $d(n)$ is our $d(n+1)$.)

Divisor methods can be described in several different ways that yield the
same result. 
The methods have been invented and reinvented in different guises;
tradition varies, 
and different types of formulations are used in, for example, election laws and
other official and non-official descriptions from different countries. There
are two main types of formulation, and we give each in two different versions.
(The equivalence of the different types of formulations have been known
for a long time.
For example,
D'Hondt proposed his method in \cite{D'Hondt} using essentially our
formulation \ref{Diva} below, 
and later showed the equivalence
to \ref{Divc}, see \cite[p.\ 125]{Kopfermann}.)

In the first type of formulation, 
the method is seen as a way of rounding.
Given a sequence $d(1), d(2),\dots$ as above, we define the
\emph{$d$-rounding} of a real number $x>0$ by
\begin{equation}
  \label{dround}
\rd{x}\=n \qquad \text{if} \quad d(n)\le x\le d(n+1),
\end{equation}
where $d(0)=0$.
Note that this is unambiguous only if  $d(n)<x<d(n+1)$; if $x=d(n)$ for an
integer $n$, then both $\rd x=n$ and $\rd x=n-1$ are acceptable values.
(This is important in the case of ties, see below.)
In particular, 
$d(n)=n-\frac12$ (Webster's method) yields standard rounding; 
$d(n)=n$ (Jefferson's method) yields rounding downwards;
$d(n)=n-1$ (Adams's method) yields rounding upwards.
(Other choices of $d(n)$ yield non-standard roundings. Note that
``rounding'' in general should be interpreted in a very weak sense,
and that $|\rd x-x|>1$ may be possible, 
see \refR{Rexact}.)
In this formulation, the numbers $d(n)$ are sometimes called
\emph{signposts},
see \cite{BY} and \cite{Pukel2008b}.

\begin{table}
  \begin{tabular}{l|l l}
& $d(n)$ & decimal values (rounded)\\
\hline
Jefferson, D'Hondt & $n$ & 1, 2, 3, 4,  \dots \\[2pt]	
Webster, \StL & $n-\frac12$ & 
0.5, 1.5, 2.5, 3.5, \dots \\[2pt]	
Adams & $n-1$ & 0, 1, 2, 3,  \dots \\[2pt]	
Imperiali & $n+1$ & 2, 3, 4, 5,  \dots \\[2pt]
Danish & $n-\frac23$ & 
0.333, 1.333, 2.333, 3.333, \dots \\[2pt]	
Adjusted \StL &
$\begin{cases}
0.7, &  n=1 \\[-1pt]	
n-\frac12,   & n\ge1
\end{cases}$
& 0.7, 1.5, 2.5, 3.5, \dots
\\[2pt]	
Cambridge Compromise & $(n-6)_+$ & 0, 0, 0, 0, 0, 0, 1, 2, 3,  \dots \\[2pt]	
Huntington & $\sqrt{n(n-1)}$ & 0, 1.414, 2.449, 3.464, \dots \\[2pt]	
Dean & $\frac{2n(n-1)}{2n-1}$  & 0, 1.333, 2.400, 3.429, \dots \\[2pt]
Estonia & $n^{0.9}$ & 1, 1.866, 2.688, 3.482, \dots \\[2pt]	
Macau & $2^{n-1}$ & 1, 2, 4, 8, \dots 
  \end{tabular}
\caption{Some divisor methods. (For Dean's method, $d(n)$ is the harmonic
  mean of $n$ and $n-1$; \cf{} Webster's and Huntington's methods with the
  arithmetic and geometric means.)}
\label{TabDiv}
\end{table}

Using the concept of $d$-rounding, and the notation from \refS{Snotation},
the divisor method may be defined as follows:

\begin{divisor}
  \label{Diva}
Let
\begin{equation}\label{diva}
  s_i\=\rd{\frac{v_i}{D}},
\end{equation}
where $D$ is a (real) number chosen such that $\sumim s_i=N$.
\end{divisor}

The rationale for this formulation is that $D$ is regarded as the price of a
seat, \ie{} the number of votes that a seat ``costs'', which would give $v_i/D$
seats to party $i$; this real number is $d$-rounded to an integer to obtain
$s_i$. 
The price $D$ is set by ``the market'' (\ie{} the election officer, or
computer) so that 
the desired total number of seats is distributed.
This interpretation can be further combined with the definition
\eqref{dround} to say that the price of $n$ seats is $d(n)D$ votes.

\begin{remark}\label{Rd(1)}
In particular, \eqref{dround} yields $\rd x=0$ if and only if $x<d(1)$,
or perhaps $x=d(1)$.
Thus, if $d(1)=0$, then $\rd x\ge 1$ for every $x>0$, while if $d(1)>0$, then
$\rd x=0$ for small positive $x$. 

Hence,  if $d(1)=0$, then \eqref{diva} yields
$s_i\ge1$, so every party (with at least
  one vote) is guaranteed at least one seat.
This is unacceptable in general elections, where on the contrary
there often are special threshold rules to prevent small parties from
getting seats. However, this is acceptable, and may even be desirable, when
distributing seats 
between states or constituencies, where typically there is a side condition of
at least one seat each (and sometimes a minimum of two or more seats); 
a typical example is the
United States Constitution \cite{BY}. 
(If $d(1)=0$, we assume that $N\ge m$, so that there are enough seats to
satisfy this minimum requirement.)

A divisor method is by \cite{Pukel2008b} 
called \emph{impervious} if $d(1)=0$ and \emph{pervious}
if $d(1)>0$. 
\end{remark}

If we move $D$ continuously from $\infty$ down to 0, it is seen that the
total number of seats $\sumim s_i$ given by \eqref{diva}
grows monotonously from 0 (when $d(1)>0$)
or $m$ 
(when $d(1)=0$) to $\infty$, and there is a value of $D$ that yields the
desired sum $\sumim s_i=N$. In general there is an interval of such $D$, all
yielding the same seat distribution $s_i$, and this determines all $s_i$
uniquely. In exceptional cases, there is only a single $D$ that works; in
this case we must have equalities $v_i/D=d(s_i)$ for some party $i$
and $v_j/D=d(s_j-1)$ for some other party $j$, see \eqref{diva} and
\eqref{dround}, and in this case the result $(s_i)_1^m$ is not unique; we
may move one seat from party $j$ to party $i$ and \eqref{diva} still holds. 
Such cases are known as \emph{ties}; they are usually resolved by drawing
lots, but other special rules for ties are used in some countries.
The formulation \ref{Diva} thus gives a well-defined method in the sense
that the numbers $s_i$ are uniquely determined, except possibly when there
are ties. (Note that the number $D$ is not uniquely determined.)

It is sometimes, as in \refS{SpfD}, convenient to rewrite \eqref{diva} as
follows using $x\=1/D$ or $t\=V/D$.
A divisor method may thus also be called a
\emph{multiplier method}, as in \cite{HPS2005}.
\begin{divisor}
  \label{Divam}
Let
\begin{equation}\label{divam}
  s_i\=\rd{{v_i}{x}}=\rd{{p_i}{t}},
\end{equation}
where $x$ or $t$ is a (real) number chosen such that $\sumim s_i=N$.
\end{divisor}

Formulations of the type \ref{Diva} of divisor methods are the standard 
in USA, where Thomas Jefferson formulated
a method of this type (with rounding downwards) in 1792
for apportionment; Jefferson's method was
adopted by Congress in 1792  and used until 1832 \cite{BY}.
Four other divisor methods 
(using \refform{Diva} with different $d(n)$, see \refTab{TabDiv})
have also been important in USA;
Huntington's method is used since 1941 and
Webster's method has been used earlier,
while Adams's and  Dean's methods have never been used but have frequently
figured in discussions, see  \cite{Huntington}
and \cite{BY}.

In Europe, this type of formulation is unusual
(although \cite{CambridgeCompromise} and the current German election law
\cite[\S{} 6 (2)]{DEvallag}
are examples), and divisor methods are more commonly
defined by formulations such as the two following (closely related) ones.
Note that these formulations are algorithmic, in the sense that they
directly yield workable methods, while with \ref{Diva} one
has to search for a suitable divisor. (Preferably by computer, where such
searches are easy.)
The European tradition of divisor methods is younger than the American one,
see \eg{} \cite{Carstairs}, \cite{Farrell}, \cite{Kopfermann}.
D'Hondt's method (which is the same as Jefferson's) was
proposed by D'Hondt \cite{D'Hondt} in 1882 and first adopted in Belgium
1899; it is now used in many countries.
\StL's method (which is the same as Webster's)
was proposed by \StL{} \cite{StL1,StL2} in 1910, and is also
used in several places.


\begin{divisor}
  \label{Divb}
The $N$ seats are awarded sequentially, with each seat given to the party $i$
that has the highest value of the quotient $v_i/d(s_i+1)$ where $s_i$ is the
number of seats that the party has received so far.
\end{divisor}

The quotient
$v_i/d(s_i+1)$ is called the \emph{comparative figure} of the party. Note that
the comparative figure is updated (decreased) each time the party receives a
seat.

\begin{divisor}
  \label{Divc}
Divide each $v_i$ by the numbers $d(1)$, $d(2)$, \dots (as far as necessary).
Assign the seats to the $N$ largest quotients $v_i/d(j)$.
\end{divisor}

With \refform{Divc},
the quotients $v_i/d(j)$ can be arranged in a table (matrix);
alternatively, the sequences of quotients obtained for different parties can
be merged 
into a single sequence in decreasing order.
It is immediate that formulations \ref{Divb} and \ref{Divc} are equivalent;
\ref{Divb} picks the largest quotients in \ref{Divc} in decreasing order.

Of course, ties can occur in formulations \ref{Divb} and \ref{Divc} as well,
and again they are resolved by drawing lots, or possibly by some other rule.
(In \refform{Divb}, ties may also occur at intermediate stages, but they do
not affect the final result.)

To see that the formulations \ref{Diva}--\ref{Divc} are equivalent and
yield the same method, note first that \eqref{diva} can be written, using
the definition \eqref{dround},
\begin{equation}
  \label{dx1}
d(s_i)\le\frac{v_i}D\le d(s_i+1),
\qquad i=1,\dots,N,
\end{equation}
or, equivalently (interpreting $v_i/0=+\infty$)
\begin{equation}
  \label{dx2}
\frac{v_i}{d(s_i+1)}\le
D\le
\frac{v_i}{d(s_i)},
\qquad i=1,\dots,N.
\end{equation}
Consequently, $(s_i)_{i=1}^N$ are such that
\begin{equation}
  \label{dx3}
\max_i \frac{v_i}{d(s_i+1)}\le
\min_i\frac{v_i}{d(s_i)};
\end{equation}
conversely, for any such $s_i$ 
we may choose $D$ such that \eqref{dx2} holds. 
Thus \refform{Diva} is equivalent to:
\begin{divisor}
  \label{Divd}
  $s_1,\dots,s_m$ are chosen such that \newline 
$\sumim s_i=N$ and \eqref{dx3} holds.
\end{divisor}
Furthermore, it is easily seen that \refform{Divd} is equivalent to \ref{Divc},
and thus also to \ref{Divb}. Hence all formulations are equivalent. 
(It is easily verified that also ties appear simultaneously in the different
formulations, except that in \refform{Divb} there may be additional
intermediate ties that do not affect the result.)
Further equivalent formulations are given by
\cite{Pukel2008b}.

We see also, from \eqref{dx2}, 
that the possible choices of $D$ in \refform{Diva} that yield
$s_1,\dots,s_m$ are $D\in[D_-,D_+]$, where
\begin{equation}
  D_-\=\max_i \frac{v_i}{d(s_i+1)}, \qquad
D_+\=\min_i\frac{v_i}{d(s_i)}.
\end{equation}
Thus, $D_+$ equals the winning comparative figure when the last seat is
awarded in formulation \ref{Divb},
and $D_-$ equals the winning comparative figure for the next seat (if there were
one).

\begin{remark}\label{RDprop}
  It follows immediately from any of the
  formulations \ref{Diva}--\ref{Divd} that the sequences $d(n)$ and
  $cd(n)$, for any constant $c>0$, define the same divisor method.
For example, \StL's method is given in \refTab{TabDiv} with $d(n)=n-\frac12$;
it can as well be defined with $d(n)=2n-1$. (In fact, the method is
traditionally 
defined with this sequence, see \eg{} \citet{StL1,StL2}.
It is therefore also called the \emph{odd-number method}.)
Similarly, the adjusted \StL's method 
(used in Sweden and Norway)
is traditionally
defined with the sequence 1.4, 3, 5, 7, \dots, and the Danish method 
with
the sequence 1, 4, 7, 10, \dots. 
(This method is used in Danish parliamentary elections for
part of the distribution of adjustment seats between consituencies;
the main distribution among parties is by the method of largest remainder.)
It may be convenient to use integers; on the other hand, 
Imperiali's method (used in Belgian local elections)
is in the election law \cite[Art. 56]{Belgien-lokal}
described using
1, 1$\frac12$, 2, 2$\frac12$, \dots, \ie, $(n+1)/2$.
In \refTab{TabDiv}, we have normalised $d(n)$ to the form $n+\gb$ whenever
possible, cf.\ \refS{ADlinear}.
\end{remark}

\begin{remark}\label{Rexact}
  Several papers impose the condition that $n-1\le d(n)  \le n$, but this is
  really not necessary.
This condition
is natural for the interpretation of the method as a kind of rounding in
\refform{Diva}, 
since it is equivalent to $\rd n=n$ for any integer $n\in\bbN$, see
\eqref{dround}. 
In particular, this means that if the exact proportions $q_i=p_iN$ happen to be
integers, then the divisor method gives $s_i=q_i$. (Take $D=V/N$ in
\eqref{diva}; \cf{} \eqref{q_i}.)
However, the method is well-defined also for other sequences
$d(n)$, and as seen in \refTab{TabDiv}, there are
some  divisor methods currently in use for general elections that do not
satisfy this condition 
(the Imperiali method \cite[Art. 56]{Belgien-lokal}
and the methods in Estonia 
\cite[\S{} 62 (5)]{EEvallag} and
Macau \cite[Artigo 17]{Macau}). 
\refForm{Diva} remains formally valid in these cases too, but in
practice these 
methods are described using \eg{} \refform{Divb} where there are no problems of
interpretation for any increasing sequence $d(n)$.
(Methods satisfying $s_i=q_i$ when all $q_i$ are integers are called
\emph{weakly proportional}
by \cite{BY}.)
\end{remark}

\begin{remark}\label{RD<=}
We required above that the sequence $d(n)$ is strictly increasing.
In fact, it can be defined assuming only $0\le d(1)\le d(2)\le\dots$, with
only a few minor modifications above.
This extension is hardly used in practice, 
except that a minimum requirement of at least $r$ seats for each party 
for some $r>1$
can be achieved  by taking $d(1)=d(2)=\dots=d(r)=0$,
\cf{} \refR{Rd(1)} for the case $r=1$. 
(We then need $N\ge rm$.)

For example, the Cambridge Compromise \cite{CambridgeCompromise}
proposed for apportionment of the European Parliament
can be described as giving each state 
a ``base'' of 
5 seats and distributing the rest by
Adams's method (subject to a maximum of 96 seats for each state); 
this is (\eg{} by formulation \ref{Diva}, see also \refT{Tgb+1} below) 
equivalent to giving each state
6 seats and distributing the rest by D'Hondt's method (subject to the same
maximum). 
If we ignore the restriction to at most 96 seats
(as done in \refTab{TabDiv} for simplicity), 
this is the same as a divisor method with $d(n)=0$ for $n\le6$ and
$d(n)=n-6$ for $n>6$, \ie{} $d(n)=(n-6)_+$. 
The maximum restriction is easily implemented if \eg{} formulation \ref{Divb}
is used, \cf{} \cite{CambridgeCompromise}; formally we can also implement it
by redefining $d(n)=\infty$ for $n>96$.
\end{remark}

\begin{remark}
  \label{Rdivisor?}
The name ``divisor method'' is used for methods of this type in any of the
formulations above. However, the word \emph{divisor} is usually used for the
(variable) number $D$ in formulation \ref{Diva}, but for the (fixed) numbers
$d(1), d(2), \dots$ in formulations \ref{Divb}--\ref{Divc}.
\end{remark}

\subsection{Linear divisor methods}\label{ADlinear}

We say that a divisor method is \emph{linear} if $d(n)=an+b$ for some $a>0$
and $b\in\bbR$. By \refR{RDprop}, we may replace $d(n)$ by $d(n)/a=n+b/a$.
We define $\gb\=b/a+1$. In other words, we may assume that $a=1$ and
\begin{equation}\label{dlin}
  d(n)=n-1+\gb
\end{equation}
for some real $\gb$. 
The reason for our choice of $\beta$ as parameter, and thus formula
\eqref{dlin},  is \eqref{dgb0} below.
Note that  $\gb=d(1)$.

We call the divisor method with $d(n)$ given by \eqref{dlin} the \emph{\gbdiv}.
The parameter $\gb$ is called \emph{proportionality index}
(\emph{Pro\-portionali\-tätsindex}) by
\cite{Kopfermann}.
The method is called \emph{$q$-stationary multiplier method} (where
$q=\gb$) by \cite{HPS2005}.

With a  linear divisor method \eqref{dlin}, the definition \eqref{dround} of
$d$-rounding yields that if $\rd x=n$, then
$n-1+\gb \le x \le n+\gb$, or equivalently
$x-\gb \le n \le x-\gb+1$. 
Consequently, $d$-rounding equals the $\gb$-rounding defined in \eqref{rga}
and \eqref{rga1}; \ie{}
\begin{equation}\label{dgb0}
  \rd x = \rgb x,
\end{equation}
at least provided $x>0$ and $x>\gb-1$.
Hence, the \gbdiv{} can, using formulation \ref{Diva}, 
be defined by \eqref{gbdiv} as we did in \refS{Snotation}.
If $\gb>1$, we here have to assume that $N$ is so large that $v_i/D>\gb-1$
and thus $\rgb{v_i/D}\ge0$. 

\begin{remark}\label{Rgb>1}
If $\gb>1$ and $N$ is small, \eqref{gbdiv}
may give a negative seat number $s_i$ for a very small party. Of course,
this has to be replaced by 0, as given by \eqref{diva}.
Hence, we may still use \eqref{gbdiv} if we ignore any party that
would get a negative number of seats. 

If $\rgb{v_i/D}\le-1$ for some $i$, then
$v_i/D\le\gb-1$ by \eqref{rga1} and thus for any party $j$, using \eqref{rga},
\begin{equation*}
 \rgb{\frac{v_j}{D}}
\le\rgb{\frac{v_j}{v_i}(\gb-1)} 
=\rgb{\frac{p_j}{p_i}(\gb-1)} 
\le \Bigpar{\frac{p_j}{p_i}-1}(\gb-1), 
\end{equation*}
and summing over $j$ we get $N\le (1/p_i-m)(\gb-1)$.
Hence, if $N>(1/\min_i p_i-m)(\gb-1)$, then \eqref{gbdiv} holds without
exception.   
\end{remark}

\begin{remark}\label{Rgb<0}
The definition above requires $\gb=d(1)\ge0$.
We can extend it to $\gb<0$ by replacing \eqref{dlin} by 
\begin{equation}\label{dlin-}
d(n)=(n-1+\gb)_+,  
\end{equation}
\ie, letting $d(n)=0$ for $n\le 1+\floor{\abs{\gb}}$, \cf{} \refR{RD<=}.
Then \eqref{dgb0} still holds for $x>0$, and \eqref{gbdiv} holds as before.
By \refR{RD<=}, every party now gets at least $\floor{|\gb|}+1$ seats, and
thus we need $N\ge m\bigpar{\floor{|\gb|+1}}$.

The Cambridge Compromise (without maximum rule) is an example, with $\gb=-5$.
\end{remark}

\begin{remark}\label{Rgb01}
Linear divisor methods are often considered only for the case $0\le\gb\le1$,
which is equivalent to the condition $n-1\le d(n)\le n$ discussed in
\refR{Rexact}. However, we see that they are
well-defined for arbitrary real $\gb$, although some care may be needed for
small $N$. 
Note that at least the case $\gb=2$ is used in practice.
\end{remark}

In the main part of the paper we consider only linear divisor methods.
Note that many of the methods in \refTab{TabDiv} are linear:
Jefferson/D'Hondt ($\gb=1$); 
Webster/\StL{} ($\gb=1/2$); 
Adams ($\gb=0$); 
Imperiali ($\gb=2$); 
Danish ($\gb=1/3$);
Cambridge Compromise ($\gb=-5$).

\begin{remark}
The adjusted \StL{} method 
also has a linear $d(n)$ as in \eqref{dlin} ($\gb=1/2$) 
except for $d(1)$, which does not affect asymptotic results.
More precisely, the adjusted \StL{} method 
differs from \StL's method only in the value of $d(1)$, and
$\rd x$ is the same for both methods for all $x\ge0.7$; hence
they give exactly the same result as soon as every party gets at least one
seat by the adjusted method. (But the adjustment
makes it more difficult for a small party to get the first seat.) For
asymptotic results as in the present paper, there is thus no difference
between the adjusted \StL{} method and \StL's method.  
\end{remark}

\begin{remark}  \label{RHuntington}
Huntington's and Dean's methods are asymptotically linear in the sense that
$d(n)=n-1+\gb+o(1)$ 
as \ntoo, in both cases with $\gb=1/2$.
An argument similar to the proof in \refS{SpfD} shows that as \Ntoox,
the probability that these methods yield the same result as Webster's tends to
1. In particular, \refT{TD} (with $\gb=1/2$) holds for these methods too, and
\refT{TDmean} shows that they are asymptotically unbiased.

Note that our asymptotic approach does not answer the controversial question
whether Webster's or Huntington's method is the most fair and unbiased, see
\cite{BY}; this depends mainly on how one measures the bias for small parties
(states), in particular for the ones obtaining just 1 seat.  
\end{remark}

Although linear divisor methods may be biased, see \refS{Sresults}, they are
perfectly proportional with respect to changes in the total number of seats.
(If $\gb<0$ or $\gb>1$, we assume that $N$ is not too small, see Remarks
\ref{Rgb>1}--\ref{Rgb<0}.) 

\begin{theorem}
  \label{TDND}
Consider a \gbdiv.
If the number of seats is increased from $N$ to $N+L$, and $Lp_1,\dots,Lp_m$
all are integers, then party $i$ gets exactly $Lp_i$ seats more:
\begin{equation}
  s_i(N+L)=s_i(N)+Lp_i.
\end{equation}
\end{theorem}

\begin{proof}
	This is easiest seen using the multiplier formulation \ref{Divam}.
If $t$ yields $s_i(N)=\rgb{p_it}$ with $\sumim s_i(N)=N$, then 
$t+L$ yields 
\begin{equation*}
s_i=\rgb{p_i(t+L)}=\rgb{p_it}+p_iL=s_i(N)+p_iL,
\end{equation*}
with $\sumim s_i=\sumim s_i(N)+\sumim p_iL=N+L$.
\end{proof}

There is also a simple relation between the methods for two different values
of $\gb$ that differ by an integer. 

\begin{theorem}\label{Tgb+1}
  Consider the \gbdiv{} and suppose that $N$ is so large that every party
  gets at least one seat. (Otherwise, ignore the remaining parties.)
Then the \gbdiv{} yields the same result as first giving one seat to each
party and then distributiong the remaining $N-m$ seats by the
$(\gb+1)$-linear divisor method.
\end{theorem}

\begin{proof}
  This follows immediately from \refform{Divc}, noting that by \eqref{dlin},
  or more generally \eqref{dlin-}, $d_{\gb+1}(n)=d_\gb(n+1)$ (using
  $\gb$ as a subscript to denote the divisor sequences).

Alternatively, this follows from \eqref{gbdiv} and the relation 
$\rgx{\gb+1}{x}=\rgx{\gb}{x}-1$, which follows from \eqref{dga1}.
\end{proof}

By repeated applications of this theorem, we can reduce any \gbdiv{}
to the case
$0<\gb\le1$. 

\begin{example}
Adams's method ($\gb=0$)
is the same as giving each party one seat and
distributing the rest by Jefferson's method $(\gb=1$).
\end{example}

\begin{example}
The Imperiali method ($\gb=2$)
is equivalent to using D'Hondt's method ($\gb=1$)
but retracting one seat from each party and redistribute them
(among the parties that had at least one seat) by continued application of
D'Hondt's method.
\end{example}

\begin{example}
The Cambridge Compromise ($\gb=-5$), is as said in \refR{RD<=}
the same as giving each party (state) 5 seats and
distributing the rest by Adams's method ($\gb=0$),
or as giving each party (state) 6 seats and
distributing the rest by 
Jefferson/D'Hondt's method $(\gb=1$).
\end{example}

\subsection{Divisor methods and proportionality}\label{SSdiv-prop}

We make a formal definition of proportionality for election methods, using
an asymptotic property, since no election method is exactly proportional
except in exceptional cases.
(\citet{BY} uses a different definition, \cf{} \refR{Rexact}.)

\begin{definition}  \label{Dprop}
  An election method is 
\emph{(asymptotically) proportional} if, for any number $m$ of
  parties and any proportions of votes $\ppm$,  
for each party $i$,
\begin{equation}\label{sin}
s_i(N)/N\to p_i
\qquad\text{as \Ntoo}.
\end{equation}
\end{definition}

Many proportional methods that are used are divisor methods, but not all
divisor methods are proportional.
In order to characterise the proportional divisor methods, recall that a
positive measurable function $f$ on some interval
$(A,\infty)$ is \emph{regularly varying with
  index $\rho$}, where $\rho$ is a real number, if
\begin{equation}\label{reg}
f(\gl x)/f(x)\to \gl^\rho \qquad\text{as \xtoo},
\end{equation}
for every $\gl>0$.  In the special case $\rho=0$, \ie{} when $f(\gl
x)/f(x)\to1$ as \xtoo{} for every $\gl>0$, we say that $f$ is \emph{slowly
  varying}. 
A typical example of a slowly varying function is $\log x$. 
A sequence $c_n$ is said to be regularly varying with index $\rho$
(and slowly varying if $\rho=0$)
if the function $c_{\floor x}$ is; this is equivalent to
$c_{\floor{\gl n}}/c_n\to\gl^\rho$ as \ntoo, for every $\gl>0$.
(See \cite{BinghamGoldieTeugels}, where many more results are given.)

\begin{theorem}
  \label{TDprop}
For a divisor method defined by $d(1), d(2),\dots$, the following are
equivalent. 
\begin{romenumerate}
\item \label{tpp}
The method is proportional.
\item \label{tpq}
For any $m$ and $\ppm$, and all $i,j\le m$,
\begin{equation}\label{syn}
\frac{s_j(N)}{s_i(N)}\to \frac{p_j}{p_i}
\qquad\text{as \Ntoo}.
\end{equation}
\item \label{tpx1}
$x\mapsto \rd x$ is regularly varying with index $1$.
\item \label{tpx0}
$x\mapsto \rd x/x$ is slowly varying.
\item \label{tpd1}
The sequence $d(n)$ is regularly varying with index $1$.
\item \label{tpd0}
The sequence $d(n)/n$ is slowly varying.
\end{romenumerate}
\end{theorem}

The proof shows that it suffices that \eqref{sin}, or \ref{tpq}, holds when
$m=2$. 

\begin{proof}
\ref{tpp}$\iff$\ref{tpq}.
If \eqref{sin} holds for each $i$, then also \eqref{syn} holds.
On the other hand, if \eqref{syn} holds, then summing \eqref{syn} over all
$j$ we obtain $N/s_i(N)\to1/p_i$, and thus \eqref{sin}.  

\ref{tpq}$\iff$\ref{tpx1}.
By \eqref{diva}, 
$s_i=\rd{v_i/D}=\rd{p_ix}$ with $x\=V/D$, and hence
\eqref{syn} can be written
\begin{equation}\label{synx}
\frac{\rd{p_jx}}{\rd{p_ix}}\to \frac{p_j}{p_i}
\qquad\text{as \xtoo},
\end{equation}
since $N=\sumim s_i=\sumim \rd{p_ix}$ and thus $N\to\infty\iff \xtoo$. 
If $\rd x$ is regularly varying with index 1, then \eqref{synx} follows from
\eqref{reg} with $\gl=p_j/p_i$. Conversely, if 
\ref{tpq}
holds, and $\gl>0$, consider the case $m=2$ and $p_1=\gl/(1+\gl)$,
$p_2=1/(1+\gl)$. Then \eqref{synx} with $(j,i)=(1,2)$ and $y\=p_2x$
yields $\rd{\gl y}/\rd y\to \gl$, which is \eqref{reg} with $\rho=1$.

\ref{tpx1}$\iff$\ref{tpx0}.
Directly by the definition \eqref{reg}.

\ref{tpx1}$\iff$\ref{tpd1}.
Let $g(x)$ be the inverse of $\rd x$ defined by
$g(x)\=\inf\set{y:\rd y>x}$. Then \eqref{dround} shows that
$g(x)=d(\floor{x}+1)$. \cite[Theorem 1.5.12]{BinghamGoldieTeugels} shows
that if $\rd x$ is regularly varying with index 1, then so is $g(x)$. This
implies $d(n+1)/d(n)\to1$ as \ntoo, and it follows that $d(n)$ is regularly
varying with index 1. 
The converse follows in the same way, now using the inverse of
$d(\floor{x})$ which is $\rd x+1$.

\ref{tpd1}$\iff$\ref{tpd0}.
Directly by the definition.
\end{proof}

In particular, \refT{TDprop} shows that any divisor method with 
\begin{equation}
  \label{prope}
d(n)/n\to a \qquad\text{as \ntoo},
\end{equation}
for some $a>0$, is proportional.
Of the methods in \refTab{TabDiv}, the only methods that do not satisfy
\eqref{prope}, and thus are proportional, are the methods of Estonia and
Macau. In fact, it is easily seen from \refform{Divb} or \ref{Divc} that the
Estonian method is the same as D'Hondt's method applied to 
$v_i^{1/0.9}=v_i^{1.111\dots}$; thus it
makes $s_i$ proportional to $v_i^{1.111\dots}$. 
This gives a small but clear
advantage to larger parties.
Macau's method, on the contrary, favours small parties and encourages
splintering. (With only 12 directly elected members of the Legislative
Assembly of Macau, it essentially imposes a maximum of 2 seats per party.)

In fact, the proportional divisor methods in \refTab{TabDiv}, and all
linear divisor methods, satisfy the  following stronger form of proportionality.

\begin{theorem}\label{TDprop2}
  For a divisor method with $d(n)=an+O(1)$ for some $a>0$,
$s_i(N)=p_iN+O(1)$. (The implicit constant may depend on $m$, but not on
  anything else.)
\end{theorem}

\begin{proof}
By \refR{RDprop}, we may replace $d(n)$ by $d(n)/a$, and thus assume
$d(n)=n+O(1)$. 

  If $d(n)=n+O(1)$, 
then also $d(n-1)=n+O(1)$ and thus \eqref{dround} shows that if $\rd x=n$,
then $x=n+O(1)$; in other words,
\begin{equation}
  \rd x = x+O(1).
\end{equation}
Consequently, \eqref{diva} yields, with $y\=V/D$,
\begin{equation}\label{db1}
  s_i=\frac{v_i}D+O(1)={p_i}y+O(1).
\end{equation}
Summing over $i$ we find $N=\sumim p_i y+O(1)=y+O(1)$, and thus \eqref{db1}
yields $  s_i={p_i}N+O(1)$.
\end{proof}

\subsection{Uniformity}\label{Auniform}
Divisor methods satisfy the following consistency property, 
called \emph{uniformity} in \cite{BY},
when we consider
subsets of the set of all parties. This 
is obvious from any of the formulations \ref{Diva}--\ref{Divd}, and
would perhaps not be worth mentioning, except for the fact that it does
\emph{not} hold for quota methods;
in fact, under some weak extra conditions, it is  satisfied only by divisor
methods \cite{BY}. 

\begin{theorem}\label{Tconsistency}
Suppose that some set of parties $A_1,\dots,A_\ell$ (where $1\le\ell\le m$) get
together $N'$ seats in an election by a divisor method.
Then the number of seats for each of these parties is the same as if $N'$
seats were distributed, by the same divisor method, 
in an election with only these parties participating
(and obtaining the same numbers of votes $v_i$).
\nopf
\end{theorem}

\section{Quota methods}\label{Aquota}

In a \emph{quota method}
one first determines a \emph{quota} $Q$;
different quota methods differ in the formula for $Q$, see below.
($Q$ can be seen as the standard price of a seat, just as the divisor $D$ in
the \refform{Diva} of divisor methods.)
The seats are then distributed as follows:
\begin{kvot}\label{kvota}
Divide the numbers of votes $v_i$ by $Q$, and give first each party 
as many seats as
the integer part $\floor{v_i/Q}$ of its fraction.
Any remaining seats are given to the parties that have largest fractional part
$\frax{v_i/Q}=v_i/Q-\floor{v_i/Q}$.
(Or, equivalently, largest remainder at the division
$v_i/Q$, since the remainder is
$v_i-\floor{v_i/Q}Q$, which is  $Q$ times the fractional part
$v_i/Q-\floor{v_i/Q}$.)
\end{kvot}

The most common quota method is 
the \emph{method of largest remainder}, or
\emph{Hamilton's method},
in Europe often called 
\emph{Hare's method},
which uses the \emph{simple quota}, 
also called
\emph{Hare quota}, 
$  Q\=\xfrac{V}{N}$,
\ie{} the average number of votes per seat.
Note that then $v_i/Q=q_i$, the exactly proportional real allocation, see
\eqref{q_i}. 
(Alexander Hamilton suggested the method in 1792 for apportionment to the US
Congress; 
it was approved by Congress but vetoed by president Washington; it made a
comeback and was used (under the name \emph{Vinton's method})
1850--1900, see \cite{BY}. 
The namn Hare's method is also well-established but is a misnomer;
Thomas Hare really advocated a different method, the 
\emph{Single Transferable Vote, STV},
which is not based on party lists and is not
treated here.) 

Another version is
\emph{Droop's method} which uses the \emph{Droop quota}
$  Q\=\xfrac{V}{(N+1)}$.
Also the \emph{Imperiali quota}
$  Q\=\xfrac{V}{(N+2)}$
has occasionally been used 
(\eg{} in Italy until 1993 \cite{Farrell}).

We define, more generally, for any real number $\gam$, the \emph{\gamqm} to
be the quota method with quota 
\begin{equation}\label{gamQ}
 Q\=\frac{V}{N+\gam}.
  \end{equation}
(For $\gam<0$, we assume that $N>|\gam|$.)
Thus Hamilton/Hare's method is the case $\gam=0$,
Droop's method is $\gam=1$, and the Imperiali quota is $\gam=2$.
(\citet{Kopfermann} calls the \gamq{} method
\emph{Rundungsverfahren mit dem Proportionalitätsindex $\rho$},
where his $\rho=(\gam+1)/2$, \ie, $\gam=2\rho-1$.)

\begin{remark}
\label{Ravrundning}
The quotas above are sometimes rounded to integers (upwards, downwards or by
standard rounding); see \cite{Pukel2010d} for several different examples in
current use. There is no mathematical reason to round the quota; on the
contrary, there are reasons against. 
\citet{Kopfermann} points out that quota methods involve rounding at the end,
but that, as a general principle,
intermediate values in a calculation should not be rounded.
A specific problem is that rounding the quota
means that the method is no longer homogeneous (see
\refS{Snotation}).
Moreover, with a rounded quota, in some cases a party may lose a seat by
getting an additional vote, which ought to be unacceptable; see \cite{SJV4}
for a simple example with Droop's method (although there formulated for STV,
which in the example gives the same result).
In a general election with a large number of votes, and therefore a large
quota, rounding the quota will usually make no difference, but in the
exceptional cases when it does make a difference, it can be harmful and lead 
to undesirable results.

In the present paper, we consider the methods with unrounded quotas.
\end{remark}

\begin{remark}\label{Rkvot}
In formulation \ref{kvota}, it is implicitly assumed that
$N-m\le\sumim \floor{v_i/Q}\le N$, 
so that the number of remaining seats is non-negative and not larger than
$m$, the number of parties.
For the \gamqm, this  is always satisfied if $-1\le\gam<1$ (including
Hamilton/Hare), and also if $\gam=1$ (Droop) except in a very special case
when all parties tie for the last seat. For $\gam<-1$ or $\gam>1$, 
formulation \ref{kvota} may have to be amended; in such cases
we interpret it by using one of the equivalent
formulations \ref{kvotb}--\ref{kvotd} below.
(These formulations  are well-defined for all $Q$ and 
$N$, up to the usual possibility of ties.)
\end{remark}

The rule in \ref{kvota} that the remaining seats are distributed to the
parties with largest remainders 
$v_i/Q-\floor{v_i/Q}$ can also be expressed by saying that we find a
suitable $\ga$ and round $v_i/Q$ upwards if the fractional part is greater than
$\ga$, and downwards if the fractional part is smaller. 
In other words, 
recalling the notion of $\ga$-rounding in \refS{Snotation},
a quota method can equivalently be described as follows.

\begin{kvot}\label{kvotb}
\begin{equation}
s_i\=\rga{\frac{v_i}Q},
\end{equation}
where $\ga$ is chosen such that $\sumim s_i=N$.
\end{kvot}

For the \gamqm{} we have $Q=V/(N+\gam)$, and thus
$v_i/Q=p_iV/Q=p_i(N+\gam)$.
Hence we have also the following formula.

\begin{kvot}\label{kvotc}
The \gamqm{} is given by
\begin{equation}\label{kvotcx}
s_i\=\rga{(N+\gam)p_i},
\end{equation}
where $\ga$ is chosen such that $\sumim s_i=N$.
\end{kvot}

Formulations \ref{kvotb}--\ref{kvotc} are the same as \eqref{gamq}.

Quota methods can also be described using a sequential allocation of seats
by comparative figures as in formulation \ref{Divb} of divisor methods.
(This is not the usual way to define quota methods, but it appears in 
proofs in \cite{StL2} and \cite{Polya1918}.)
It is easy to see that \ref{kvotb} is equivalent to the following.

\begin{kvot}  \label{kvotd}
The $N$ seats are awarded sequentially, 
with each seat given to the party $i$
that has the highest value of the difference $v_i-s_iQ$ where $s_i$ is the
number of seats that the party has received so far.
\end{kvot}

It is easy to see that every \gamqm{} satisfies $s_i(N)=Np_i+O(1)$, uniformly
in all $m$ and $\ppm$, and hence is proportional, \cf{} \refD{Dprop} and
\refT{TDprop2}. We omit the easy proof since we show a precise version of
this in \refT{TQbound}.

Moreover, a \gamqm{} is
perfectly proportional with respect to changes in the total number of seats,
just as a \gbdiv, \cf{} \refT{TDND}.

\begin{theorem}
  \label{TDNQ}
Consider a \gamqm.
If the number of seats is increased from $N$ to $N+L$, and $Lp_1,\dots,Lp_m$
all are integers, then party $i$ gets exactly $Lp_i$ seats more:
\begin{equation}
  s_i(N+L)=s_i(N)+Lp_i.
\end{equation}
\end{theorem}

\begin{proof}
Obvious from \eqref{kvotcx}, since the same $\ga$ works for $N$ and $N+L$.
\end{proof}

\subsection{Two parties}\label{SS2}
In the simple case $m=2$, when there are only two parties, the \gbdiv{} and
the \gamqm{} with $\gam=2\gb-1$ coincide \cite[Satz 6.2.6]{Kopfermann}. 
(In particular, then Webster/\StL{} = Hamilton/Hare and Jefferson/D'Hondt =
Droop.) 
In fact, it is
easy to see that
\begin{equation}\label{2party}
  s_i\=\rgb{(N+2\gb-1)p_i}, \qquad i=1,2,
\end{equation}
gives $s_1+s_2=N$, (at least for some choice in case of a tie).
Thus, for the \gbdiv{} we can take $D'=1/(N+2\gb-1)=1/(N+\gam)$ and 
$D=D'V=V/(N+2\gb-1)$ in \eqref{gbdiv}, and for the \gamqm{} we can take
$\ga=\gb$ in \eqref{gamq}. Hence, both methods yield the result given by
\eqref{2party}. 

\section{Lemmas on uniform distribution}\label{Auni}

\begin{proof}[Proof of \refL{Lux}]
First, by \refL{LW}, the sequence
\begin{equation}
(\frax{ny_1/m+a_1/m},  \dots,\frax{ny_k/m+a_k/m})
\end{equation}  
is uniformly distributed in $[0,1)^k$.
Thus, by the definition of $\Modm$,
\begin{equation}\label{qsw}
\bigpar{\Modm\xpar{ny_1+a_1},  \dots,\Modm\xpar{ny_k+a_k}}
\end{equation} 
is uniformly distributed in $[0,m)^k$. 
If a sequence $(x_n)$ is uniformly distributed in $[0,m)$, then 
the  sequence of pairs $(\frax{x_n},\floor{x_n})$ 
is uniformly distributed in $[0,1)\times\set{0,\dots,m-1}$. 
This extends to $k$ dimensions, and thus it
  follows from \eqref{qsw} that if
  $\tell_{nj}\=\Modm(\floor{ny_j+a_j})=\floor{\Modm(ny_j+a_j)}$,
then 
\begin{equation}\label{qm}
\bigpar{\frax{ny_1+a_1},  \dots,\frax{ny_k+a_k},
\tell_{n1},\dots,\tell_{nk}}
\end{equation} 
is uniformly distributed in $[0,1)^k\times\set{0,\dots,m-1}^k$. 

The argument just given shows
that \eqref{qm} is uniformly distributed also along any subsequence
$n=m\nu+n_0$, $\nu\ge1$.
We have $\ell_n=\Modm(n-\sumjk \tell_{nj})$, and it follows that 
along any such subsequence, $\ell_n$ is determined by
$\tell_{n1},\dots,\tell_{nk}$. Consequently,
the vectors \eqref{qulu} are
uniformly distributed along any such subsequence, and therefore also along
the full sequence.
\end{proof}

The following lemma is essentially taken from \cite{HPS2004}. 
(There Riemann integrability is assumed, but as the proof
below shows, Lebesgue integrability suffices.)

\begin{lemma}
  \label{Luni}
Let $X=(X_1,\dots,X_k)$ be a random variable with an absolutely continuous
distribution in $\bbR^k$.
Let $\nu_n$ be any sequence of constants with $\nu_n\to\infty$. Then
\begin{equation}\label{luni}
  \bigpar{\frax{\nu_nX_1},\dots,\frax{\nu_nX_k},X_1,\dots,X_k}
\dto
  \bigpar{U_1,\dots,U_k,X_1,\dots,X_k}
\end{equation}
where $U_1,\dots,U_k\sim\U(0,1)$ are independent of each other and of $X$.
\end{lemma}

\begin{proof}
The assumption that $X$ is absolutely continuous means that it has a density
function $f$, which is a Lebesgue integrable function on $\bbR^k$.

  We regard $\frax{\nu_n X_j}$ and  $U_j$ as elements of the circle group
  $\bbT=\bbR/\bbZ$. Thus the random vectors in \eqref{luni} are elements of
the group $\bbT^k\times\bbR^k$, and by standard Fourier analysis, it
suffices to show that for any integers $\ell_1,\dots,\ell_k$ and real
numbers $t_1,\dots,t_k$,
\begin{equation}\label{rla}
  \E e^{\sum_j 2\pi\ii\ell_j\frax{\nu_nX_j}+\sum_j\ii t_jX_j}
\to
  \E e^{\sum_j 2\pi\ii\ell_jU_j+\sum_j\ii t_jX_j}.
\end{equation}
Since $\nu_nX_j-\frax{\nu_nX_j}$ is an integer, we have, with $f$ as above,
\begin{multline}\label{rlb}
  \E e^{\sum_j 2\pi\ii\ell_j\frax{\nu_nX_j}+\sum_j\ii t_jX_j}
=
  \E e^{\sum_j 2\pi\ii\ell_j\nu_nX_j+\sum_j\ii t_jX_j}
\\
=
  \E e^{\sum_j \ii(2\pi\ell_j\nu_n+t_j)X_j}
=\widehat f(2\pi\ell_1\nu_n+t_1,\dots,2\pi\ell_k\nu_n+t_k).  
\end{multline}
If $(\ell_1,\dots,\ell_k)=(0,\dots,0)$, then \eqref{rla} is trivial.
If $(\ell_1,\dots,\ell_k)\neq(0,\dots,0)$, then 
$|(2\pi\ell_1\nu_n+t_1,\dots,2\pi\ell_k\nu_n+t_k)|\to\infty$, and thus 
\eqref{rlb} and the Riemann--Lebesgue lemma show that the left-hand side of
\eqref{rla} tends to 0, which verifies \eqref{rla} in this case too,
since 
$$\E e^{\sum_j 2\pi\ii\ell_jU_j+\sum_j\ii t_jX_j}
=
\prod_j \E e^{2\pi\ii\ell_j U_j}\cdot\E e^{\sum_j\ii t_jX_j}
=0
$$
because $\E e^{2\pi\ii\ell_jU_j}=0$ when $\ell_j\neq0$.
\end{proof}

\section{Moments of order statistics for random $p_i$}\label{Apox}

Let $\ppm$ be random and uniformly distributed on $\sxm$ as in
\refS{Srandom}, and let $p\ox1\le \dots\le p\ox m$ be the order statistics.
We give here, for completeness, a calculation of  moments of
$p\ox k$ by a well-known method.
Note that we here follow standard convention for order statistics
in probability theory and order $p_i$
in increasing order, in contrast to
\refS{Srandom} where we order them in decreasing order; 
thus $p\oy k=p\ox{m+1-k}$, which should be remembered when using the results
below in \refS{Srandom}.

Let $T_1,\dots,T_m$ be $m$ independent identically distributed random
variables with the exponential distribution $\Exp(1)$, and let $Z\=\sumim T_i$.
Then $(T_1/Z,\dots,T_m/Z)$ is uniformly distributed on the simplex $\sxm$, and
moreover independent of $Z$. This means that we can take $p_i=T_i/Z$. For
the corresponding order statistics we thus have $p\ox k=T\ox k/Z$; moreover
$p\ox k$ and $Z$ are independent and thus, for any $k$ and $\ell$,
\begin{equation}
  \E (T\ox k)^\ell =   \E (Zp\ox k)^\ell =\E(Z^\ell)\E(p\ox k^\ell),
\end{equation}
so the moments of $p\ox k$ are given by
\begin{equation}\label{kva2}
  \E (p\ox k)^\ell=
\frac{ \E  (T\ox k)^\ell}{ \E Z^\ell}.
\end{equation}
Furthermore, $Z\=\sumim T_i$ has the Gamma distribution $\gG(m)$ with moments
\begin{equation}\label{kva3}
  \E Z^\ell = \frac{\gG(m+\ell)}{\gG(m)}=m(m+1)\dotsm(m+\ell-1).
\end{equation}

To find the moments of the order statistics $T\ox k$, we regard
$T_1,\dots,T_m$ as the times of the first events in $m$ independent Poisson
processes. (Alternatively, we may regard them as life-lengths of $m$
identical radioactive atoms.) It is then well-known, and  easy to see by the
lack of memory in exponential distributions, that the 
increments (or waiting times) $V_j\=T\ox j-T\ox{j-1}$ (with $T\ox0\=0$) are
independent exponential random variables with $V_j\sim\Exp(1/(m+1-j))$.
Hence $T\ox k=\sumjk V_j$ with these $V_j$, and moments are easily computed.
In particular,
\begin{equation}
  \E T\ox k=\sumjk\E V_j =\sumjk\frac{1}{m+1-j}=\sum_{i=m-k+1}^m\frac1i.
\end{equation}
(This can be written as $H_{m}-H_{m-k}$, where $H_m=\sumim 1/i$ is the
$m$:th harmonic number.)
Similarly,
\begin{equation} \label{et2}
  \begin{split}
\E T\ox k^2 &=\Var T\ox k + \bigpar{\E T\ox k}^2
=\sumjk\Var V_j +\lrpar{\sumjk\E V_j}^2	
\\&
={\sumjk\frac{1}{(m+1-j)^2}}	
+
\lrpar{\sumjk\frac{1}{m+1-j}}^2	
\\&
=\sum_{i=m-k+1}^m\frac1{i^2} +\lrpar{\sum_{i=m-k+1}^m\frac1i}^2
.  \end{split}
\end{equation}
Hence, by \eqref{kva2}--\eqref{kva3},
\begin{align} \label{epox}
  \E (p\ox k) &= 
\frac{\E T\ox k }m
=\frac1m \sum_{i=m-k+1}^m\frac1i,
\\ \label{epox2}
  \E (p\ox k)^2 &= 
\frac{\E T\ox k^2 }{m(m+1)}
=\frac1{m(m+1)}\sum_{i=m-k+1}^m\frac1{i^2} 
+\frac1{m(m+1)}\lrpar{\sum_{i=m-k+1}^m\frac1i}^2 \!,
\\ \label{varpox}
  \Var (p\ox k) &= 
  \E (p\ox k)^2 -(\E p\ox k)^2 
\nonumber\\&
=\frac1{m(m+1)}\sum_{i=m-k+1}^m\frac1{i^2} 
-\frac1{m^2(m+1)}\lrpar{\sum_{i=m-k+1}^m\frac1i}^2.
\end{align}

Covariances are computed similarly.
If $1\le k\le\ell\le m$, then 
\begin{equation} 
\Cov( T\ox k,T\ox\ell)
=\sumjk\Var V_j 
=\sumjk\frac{1}{(m+1-j)^2}
=\sum_{i=m-k+1}^m\frac1{i^2}
\end{equation}
and, calculating as in \eqref{et2} and \eqref{epox},
\begin{multline}
\label{covpox}
  \Cov (p\ox k,p\ox\ell)\\
=\frac1{m(m+1)}\sum_{i=m-k+1}^m\frac1{i^2} 
-\frac1{m^2(m+1)}{\sum_{i=m-k+1}^m\frac1i}
 {\sum_{i=m-\ell+1}^m\frac1i}.
\end{multline}

\begin{example}
  For $m=3$, the covariance matrix of
  $(p\oy1,p\oy2,p\oy3)=\\
(p\ox3,p\ox2,p\ox1)$
is
\begin{equation*}
  \frac1{648}
  \begin{pmatrix}
13 & -8 & -5 \\	
-8 & \phantom{-}7 & \phantom{-}1 \\
-5 & \phantom{-}1 & \phantom{-}4
  \end{pmatrix}.
\end{equation*}
\end{example}

\newcommand\AMS{Amer. Math. Soc.}

\newcommand\vol{\textbf}
\newcommand\jour{\emph}
\newcommand\book{\emph}
\newcommand\inbook{\emph}
\def\no#1#2,{\unskip#2, no. #1,} 
\newcommand\toappear{\unskip, to appear}

\newcommand\arxiv[1]{\url{arXiv:#1.}}
\newcommand\urlq[2]{\url{#1} (#2)}

\def\nobibitem#1\par{}

\end{document}